\documentclass[12pt]{amsart}
\usepackage{graphicx}
\usepackage{amssymb}
\usepackage{amsmath}
\usepackage{amscd}
\usepackage{upgreek}
\usepackage{mathrsfs}
\usepackage{stmaryrd}
\usepackage{longtable}
\usepackage{txfonts}
\usepackage[T1]{fontenc}
\usepackage{eucal}
\usepackage{titletoc}
\usepackage{wrapfig}
\usepackage{float}
\usepackage{xypic}
\usepackage{dsfont}
\usepackage{xcolor}
\usepackage{color}
\usepackage[colorlinks,linkcolor=blue,anchorcolor=blue,
citecolor=blue]{hyperref}
\usepackage{hyperref}
\usepackage[hmarginratio=1:1,top=32mm,columnsep=20pt]{geometry}
	
\allowdisplaybreaks	
\usepackage{mathptmx}
\DeclareMathAlphabet{\mathcal}{OMS}{cmsy}{m}{n}
\DeclareSymbolFont{largesymbols}{OMX}{cmex}{m}{n}

\newtheorem{theorem}{Theorem}[section]
\newtheorem{prop}[theorem]{\rm \textsc{Proposition}}
\newtheorem{defn}[theorem]{\rm \textsc{Definition}}
\newtheorem{lem}[theorem]{\rm \textsc{Lemma}}
\newtheorem{coro}[theorem]{\rm \textsc{Corollary}}
\newtheorem{conj}[theorem]{\rm \textsc{Conjecture}}
\newtheorem{thm}[theorem]{\it \textsc{Theorem}}

\newtheorem{rem}[theorem]{\rm \textsc{Remark}}

\newcommand{\C}{\mathds{C}}
\newcommand{\id}{\textrm{id}}

\newcommand{\Z}{\mathds{Z}}

\newcommand{\g}{\mathfrak{g}}
\newcommand{\h}{\mathfrak{h}}
\newcommand{\gl}{\mathfrak{gl}}
\newcommand{\sll}{\mathfrak{sl}}
\newcommand{\psl}{\mathfrak{psl}}
\newcommand{\p}{\mathfrak{p}}
\newcommand{\gb}{\mathfrak{b}}

\newcommand{\ra}{\longrightarrow}

\newcommand{\Ind}{{\rm Ind}}

\newcommand{\ad}{{\rm ad}}
\newcommand{\Ad}{\rm Ad}

\newcommand{\B}{{\mathcal B}}

\newcommand{\A}{\mathcal{A}}
\newcommand{\aut}{{\rm Aut}}

\newcommand{\ev}{{\rm ev}}
\newcommand{\LL}{{\mathcal L}}
\newcommand{\CC}{{\mathcal C}}
\newcommand{\cZ}{{\mathcal Z}}
\newcommand{\SSS}{{\mathcal S}}
\newcommand{\OO}{\mathcal{O}}

\newcommand{\tr}{\mathfrak{tr}}
\newcommand{\ttr}{\widetilde{\mathfrak{tr}}}

\newcommand{\0}{{\bar{0}}}
\newcommand{\w}{\bar{1}}

\newcommand{\Span}{\rm span}

\begin{document}
\setlength{\oddsidemargin}{0cm}
\setlength{\evensidemargin}{0cm}

%%%%%%%%%%%%%%%%%%%%%%%%%%%%%%%%%%% TITLE, AUTHOR AND 

%%%%%%%%%%%%%%%%%%%%%%%%%%%%%%%%%%%ABSTRACT, MSC AND KEYWORDS
%\begin{abstract}
%\S 1: Sets and semigroups;
%\end{abstract}
%\subjclass[2010]{13A50.}
%\keywords{Modular invariants; complete intersection; SAGBI basis.}
%\maketitle

\title{\scshape Left-symmetric Superalgebras on 
Special Linear Lie Superalgebras}
\author{\scshape Ivan Dimitrov and Runxuan Zhang}
\address{Department of Mathematics and Statistics, 
Queen's University, Kingston K7L 3N6, Canada}
\email{dimitrov@queensu.ca}
\address{School of Mathematics and Statistics, 
Northeast Normal University, Changchun 130024, P.R. China}
\email{zhangrx728@nenu.edu.cn}
%\thanks{I. D. thanks the NSF of China (No. 11301061)  for supporting a visit to the Northeast Normal University, China and R.Z. is  grateful to the China Scholarship Council for overseas scholarships supporting her visit to Queen's University, Canada. Research was partially supported by an NSERC Discovery Grant.}}
%\date{\today}
\def\shorttitle{Left-symmetric superalgebras on $\sll(m|n)$}

%%%%%%%%%%%%%%%%%%%%%%%%%%%%%%%%%%%ABSTRACT, MSC AND KEYWORDS
\begin{abstract}
In this paper, we study the existence and classification problems of
left-symmetric superalgebras on special linear Lie superalgebras
$\sll(m|n)$ with $m\neq n$. The main three results of this paper are: 
(i) a complete classification of the left-symmetric superalgebras 
on $\sll(2|1)$, (ii) $\sll(m|1)$ does not admit 
left-symmetric superalgebras for $m \geq 3$, and (iii) 
$\sll(m+1|m)$ admits a left-symmetric superalgebra for 
every $m \geq 1$. To prove these results we combine previous results
on the existence and classification of left-symmetric algebras on 
the Lie algebras $\gl_m$ with a detailed analysis of small representations of the Lie superalgebras $\sll(m|1)$. 
We also conjecture that $\sll(m|n)$ admits left-symmetric superalgebras
if and only if $m = n+1$.
\end{abstract}


\maketitle

\baselineskip=17pt

%%%%%%%%%%%%%%%%%%%%%%%%%%%%%%%%%%%%%CONTENTS
%\textcolor{blue}{\tableofcontents{}}
\dottedcontents{section}[1.16cm]{}{1.8em}{5pt}
\dottedcontents{subsection}[2.00cm]{}{2.7em}{5pt}
\dottedcontents{subsubsection}[2.86cm]{}{3.4em}{5pt}

%%%%%%%%%%%%%%%%%%%%% let's begin %%%%%%%%%%%%%%%%%%
%%%%%%%%%%%%%%%%%%%%%%%%%%%%%%%%%%%%%%%%%%%SECTION 1
\section{Introduction} \label{secIntro}
\setcounter{equation}{0}
\renewcommand{\theequation}
{1.\arabic{equation}}
\setcounter{theorem}{0}
\renewcommand{\thetheorem}
{1.\arabic{theorem}}

A superalgebra $(\LL=\LL_{\bar 0}\oplus \LL_{\bar 1} ,\cdot)$ over a
field $k$ is called a \textit{left-symmetric superalgebra} (or an LSSA
for short) if the associator 
$(x,y,z)=(x\cdot y)\cdot z-x\cdot(y\cdot z)$ is supersymmetric in
$x$ and $y$, i.e., $(x,y,z)=(-1)^{|x||y|}(y,x,z)$; or, equivalently, 
\begin{equation}
(x\cdot y)\cdot z-x\cdot(y\cdot z)=(-1)^{|x||y|}
((y\cdot x)\cdot z-y\cdot(x\cdot z)), \;\; \forall x,y,z\in \LL.
\end{equation}
 It is clear that each associative superalgebra
is an LSSA. The supercommutator
$[x,y]=x\cdot y-(-1)^{|x||y|}y\cdot x$
defines a Lie superalgebra structure  on $\LL$. The resulting Lie
superalgebra $\mathfrak{g}_\LL$ is called the
\textit{associated Lie superalgebra} of $\LL$,  and $\LL$ is called  
an LSSA on
$\g_\LL$. Note that if $\LL$ is an LSSA on a Lie superalgebra $\g$,
then $\LL_{\0}$ is a left-symmetric algebra on $\g_{\0}$. Extending 
Segal's remark, \cite{Se1994}, we consider
determining whether the set of LSSAs on a Lie superalgebra is non-empty and
classifying all such LSSAs to be a fundamental task.

The "non-super" version of this problem arises  in the theory of affine
structures on differentiable manifolds and Lie groups. Assume
that $G$ is a connected, simply connected Lie group
with Lie algebra $\g$. It is  known that endowing $G$ with a
left-invariant  affine structure is equivalent to endowing $\g$ with
a left-symmetric product. For more details, see for example
\cite{M1981, Mi1977,  Se1994}.
Left-symmetric algebra structures on both 
finite-dimensional and  infinite-dimensional Lie algebras have been studied
extensively, see \cite{Bai2009,  Bur1994,  YZ2018} and \cite{CL2012,  KCB2011, TB2012}, respectively.
Kong and Bai, \cite{KB2008}, studied LSSAs on the Virasoro superalgebra.

In this paper we use repeatedly the following results. 
Medina, \cite{M1981}, demonstrated that finite-dimensional complex semisimple Lie
algebras  do not admit left-symmetric algebra structures.
Baues, \cite{Bau1999}, classified all left-symmetric algebras on
$\gl_n$ and proved that $\gl_n$ is the only reductive Lie algebra
with one-dimensional center and a simple semisimple ideal which admits
left-symmetric algebras over an algebraically closed ground field 
(see also \cite{Bur1996}).

Unlike the case of Lie algebras, there do exist finite-dimensional complex
simple Lie superalgebras that admit LSSAs. The problem of
classifying the LSSA-structures on Lie superalgebras is even more
challenging. Indeed, Xu in \cite{Xu2000} stated "it looks more
challenging to classify LSSAs on all the well-known simple Lie
superalgebras". Based on the classification of finite-dimensional
simple Lie superalgebras over $\C$,  the field of complex numbers, 
due to Kac,  the  even parts of classical Lie
superalgebras except for $\sll(m|n), m\neq n\geq 1$ and
$\mathfrak{osp}(2,2n), n\geq 1$ are semisimple Lie algebras. Moreover,
the even part of $\mathfrak{osp}(2,2n)$ is isomorphic to
$\mathfrak{sp}_{2n}\oplus \C$. Using the results about
left-symmetric algebras on Lie algebras, we conclude that 
$\sll(m|n), m\neq n\geq 1$ are the only classical Lie superalgebras
that may admit LSSAs.

The present paper is devoted to  the existence and  classification
of LSSAs on the Lie superalgebra $\g = \sll(m|n)$ 
with $m > n$. We completely solve these problems in the case $n =1$. 
Namely, we prove the following theorems.

\begin{thm} \label{thm1}
There do not exist LSSAs on  simple Lie superalgebras $\sll(m|1)$ 
for $m\geq3$.
\end{thm}

\begin{thm} \label{thm2}
Let $\LL$ be an LSSA on  simple Lie superalgebra $\sll(2|1)$.

\begin{enumerate}
\item $\LL$ corresponds to a bijective evaluation map associated with
an appropriate  $\sll(2|1)$-module.
\item $\LL$ is isomorphic to an LSSA in one of these three families:  
\[\A_{k}, \,\,\, k\in\C\setminus\{-1,-3\}; \quad \B_{k_1,k_2},  
\,\,\, k_1, k_2 \in \C \setminus\{0\}, k_1 + k_2 \neq -2; \quad \CC_{k},
\,\,\,  k\in\C\setminus\{0, -1\}.\]
\item $\A_k \cong \A_{-2-k}$, $\B_{k_1,k_2} \cong \B_{k_2,k_1} \cong
\B_{{-2-k_1},{-2-k_2}}\cong \B_{{-2-k_2},{-2- k_1}}$, and $\CC_k \cong
\CC_{-2-k}$.
Moreover, these are the only isomorphisms among LSSAs in (2) above.
\end{enumerate}
\end{thm}

\begin{rem} \label{rem1.10}
The families $\A_k, \B_{k_1,k_2}$, and $\CC_k$ are constructed in
Section \ref{sec5.1}.
\end{rem}

To prove these results we start with Baues' classification of left-symmetric algebras on
$\gl_n$ and then study how a left-symmetric algebra structure on $\g_{\0}$ can be
extended to
an LSSA-structure on $\g$. The basic idea is that such an
LSSA-structure 
on $\g$ exists if and only if $\g$ admits a bijective 1-cocycle
corresponding to the respective representation of $\g$.

Understanding the LSSAs on $\g=\sll(m|n)$ for $n >1$ is more difficult
because there is not a complete classification of the left-symmetric algebras on the even part 
$\g_{\0} = \sll_m \oplus \sll_n \oplus \C$. Somewhat surprisingly, we
prove the following result.

\begin{thm} \label{thm4}
There exists an LSSA on $\sll(m+1|m)$ for  every natural number $m$.
\end{thm}

Baues \cite[Proposition 5.1]{Bau1999} proved that  each left-symmetric
algebra on $\gl_n$ has a unique right identity. We conjecture that
each LSSA on $\sll(m|n)$ also has a unique right identity, which means
that $\dim \sll(m|n)_{\bar 0}=\dim \sll(m|n)_{\bar 1}$ by Proposition
\ref{prop:3.1} below. Then we have $m=n+1$. So we state the following
conjecture:

\begin{conj} \label{conj5}
There do not exist  LSSAs on any Lie superalgebras $\sll(m|n)$ other
than \hfill 
\newline $\sll(m+1|m)$. 
\end{conj}

This paper  is organized as follows. In  Section 2, we discuss the relationship between LSSAs and bijective
1-cocycles on a given Lie superalgebra. We also study evaluation maps,
which form a special class of $1$-cocycles and are useful in proving
isomorphisms of LSSAs on Lie superalgebras.
In Section 3, we present some preliminaries  on the special linear Lie superalgebras
$\sll(m|n)$. We recall the construction of Kac modules  and extensions
between irreducible modules of $\sll(m|n)$. 
Section 4 investigates $m^2|2m$-dimensional $\sll(m|1)$-modules  
for $m\geq 3$,  whose  even parts are isomorphic to the direct sum of
$m$ copies of the standard module or $m$ copies of the dual module of
the standard module as  $\sll_m$-modules.
We prove that there are no bijective $1$-cocycles of $\sll(m|1)$ for
$m\geq 3$, proving Theorem \ref{thm1}.
In Section 5, we show that bijective 1-cocycles and bijective
evaluation maps associated with $4|4$-dimensional $\sll(2|1)$-modules
coincide and classify all the LSSAs on $\sll(2|1)$. Finally, in Section
6, we prove Theorem \ref{thm4} by constructing a bijective evaluation
map on each $\sll(m+1|m)$ for $m \geq 1$. In preparation for the proof
of Theorem \ref{thm4}, we also establish some facts about representations
of (possibly infinite-diemnsional) Lie (super)algebras which may be of
independent interest.

\subsection*{{Notation and conventions}}
All vector spaces, algebras, superalgebras, etc. are over $\C$.
Elements of $\Z_2$ are denoted by $\0$ and
$\w$.  Homomorphisms (isomorphisms, automorphisms) 
of superalgebras are assumed to be homogeneous linear maps of
degree zero. 
If $W = W_{\0} \oplus W_{\w}$ is a $\Z_2$-graded vector space 
and $L : W \to W$ is a linear map, the supertrace $\mathrm{str}(L)$ of
$L$ is defined as $\mathrm{tr}(L_{\0\0}) - \mathrm{tr}(L_{\w\w})$, where
$L_{\gamma\, '\gamma\, ''} : = p_{\gamma\, ''} \circ L \circ \iota_{\gamma\, '}$
is the natural map $W_{\gamma\, '} \to W_{\gamma\, ''}$ defined by $L$, see also
Section \ref{sec3.1} for a coordinate definition of $\mathrm{str}(L)$.
A module $V$ of a superalgebra $A = A_{\0} \oplus A_{\w}$ is 
always assumed to be $\Z_2$-graded, that is $V=V_{\0}\oplus V_{\w}$ and
$A_{\gamma\, '} V_{\gamma\, ''} \subseteq V_{\gamma\, ' + \gamma\, ''}$ for
$\gamma\, ',  \gamma\, '' \in \Z_2$. We use the terms $\g$-module and representation
of $\g$ interchangeably to mean a finite-dimensional representation of
a Lie (super)algebra $\g$. 

\subsection*{{Acknowledgments}} I. D. thanks the NSFC (No. 11301061)  for supporting 
his visit to the Northeast Normal University, China and R.Z. thanks  the China Scholarship Council for
overseas scholarships for supporting her visit to Queen's University,
Canada. Research was partially supported by an NSERC Discovery Grant.
We thank the anonymous referee for their thoughtful remarks.

%%%%%%%%%%%%%%%%%%%%%%%%%%%%%%%%%%%%%%%%%%%SECTION 1
\section{Left-symmetric superalgebras} \label{sec2} 
\setcounter{equation}{0}
\renewcommand{\theequation}
{2.\arabic{equation}}
\setcounter{theorem}{0}
\renewcommand{\thetheorem}
{2.\arabic{theorem}}

 \subsection{Left regular representions} \label{sec2.1}

 Let $(\LL,\cdot)$ be an  LSSA and $(\g_{\LL}, [\;,\;])$  its
 associated Lie superalgebra. Then there are two  product operations
 $\cdot $ and $[\; ,\;]$ on the underlying $\Z_2$-graded vector space
 of $\LL$. For an element $x\in \LL$, the left multiplication operator
 $\rho(x):\LL\ra \LL$ sends  $y\in \LL$ to $x\cdot y$, and the right
 multiplication operator $\tau(x):\LL\ra \LL$ sends  $y\in \LL$ to
 $(-1)^{|x||y|}y\cdot x$. Define
$\rho:\g_{\LL}\ra \gl(\LL), \, x\mapsto \rho(x)$.
It is easy to check that $\rho([x,y])=[\rho(x),\rho(y)]$ for all
$x,y\in \g_{\LL}$, so the map $\rho$ gives a
representation of Lie superalgebra $\g_{\LL}$, which is called the
\textit{left regular representation} of $\g_{\LL}$.

The following proposition relates  LSSAs on Lie superalgebras to left
and right identities.
\begin{prop}\label{prop:3.1}
Let $\g$ be a Lie superalgebra of dimension $p|q$ satisfying
$[\g,\g]=\g$. Suppose that there exists an  LSSA  on $\g$. The following statements hold.
\begin{enumerate}
  \item $\mathrm{str} (\rho(x))=0$ and $\mathrm{str} (\tau(x))=0$ for
  all $x\in \g$.
  \item If there is an element $e\in \g$ such that
  $\rho(e)=\mathrm{id}$ or $\tau(e)=\mathrm{id}$, then $p=q$.
 \item If $\g$ is simple, then there is no element $e\in \g$ such that
 $\rho(e)=\mathrm{id}$.
\end{enumerate}
\end{prop}

\begin{proof}
(1) Since for all $x,y\in \g$, $\mathrm{str} (\rho([x,y]))=\mathrm{str}
([\rho(x), \rho(y)])=0$ and $[\g,\g]=\g$, we have
$\mathrm{str}(\rho(x))=0$ for all $x\in \g.$ Similarly,
$\tau([x,y])=\rho([x,y])-\mathrm{ad}_{[x,y]}=[\rho(x),\rho(y)]-
[\mathrm{ad}_x,\mathrm{ad}_y]$ gives $\mathrm{str}(\tau(x))=0$  for
all $x\in \g.$

(2) By (1) one has $0=\mathrm{str}(\rho(e))=\mathrm{str } (\id)=p-q$ or
$0=\mathrm{str}(\tau(e))=\mathrm{str } (\id)=p-q$, so (2) follows.

(3) Assume that there is an element $e\in \g$ such that $\rho(e)=\id$.
Since $\rho$ gives a representation of $\g$, we have that
$H=\mathrm{Ker}(\rho)$  is an ideal of $\g$. By the assumption that
$\g$ is simple  and $\rho$ is nonzero, then $H=0$. On the other hand,
$\rho([e,x])=[\rho(e), \rho(x)]=[\id, \rho(x)]=0$ for all $x\in \g$
implies $[e,\g]\subseteq H=0$, thus $e\in \cZ(\g),$ the center of $\g$.
Since $\cZ(\g)=0$, we have $e=0$ and $\rho(e)=0$, which is a
contradiction.
\end{proof}

 \subsection{LSSAs and 1-cocycles} \label{sec2.2}

Given a Lie superalgebra $\g$ and a representation $f:\g\ra \gl(V)$ of
$\g$, an even linear map $q: \g \ra V$  satisfying
\begin{equation}
\label{eq:cocy}
q([x,y])=f(x)q(y)-(-1)^{|x||y|}f(y)q(x),\;\;\forall x,y\in \g,
\end{equation}
is called an (even) \textit{1-cocycle} on $\g$  and   denoted  by the
pair $(f,q)$. A 1-cocycle $(f,q)$ is called \textit{bijective} if $q$
is a bijection.

\begin{lem}\label{lem:3.1}
Let $\g$ be a Lie superalgebra  and  $f:\g\ra \gl(V)$ a representation 
of $\g$. If $(f,q)$ is a 1-cocycle on   $\g$ such that 
$q|_{\g_{\bar 0}}=0$, then $q|_{\g_{\bar 1}}: \g_{\bar 1} \ra 
V_{\bar 1}$ is a homomorphism of $\g_{\0}$-modules.
\end{lem}

\begin{proof}
Since $q|_{\g_{\bar 0}}=0$, we have
\[q([x,y])=f(x)q(y)-f(y)q(x)=f(x)q(y), \;\; \forall x\in \g_{\bar 0},  
y\in \g_{\bar 1}.\] This means that $q|_{\g_{\w}}\circ 
\ad_x=f(x)\circ q|_{\g_{\w}}$ for all $x\in \g_{\bar 0}$, as desired.
\end{proof}

Given a Lie superalgebra $\g$, we denote by $\SSS$ the set of all LSSAs
on $\g$ and denote by $\OO$ the set of all bijective 1-cocycles on $\g$
{\footnote {Since we will be interested in bijective 1-cocycles up to 
a quasi-equivalence, see below, we may fix the vector space $V = \g$
to avoid set-theoretic pitfalls that may arise when considering 
the collection of all 1-cocycles on $\g$.}}.
Following Bai, \cite{Bai2009}, we note a close relation between $\SSS$
and $\OO$. Suppose $\SSS$ and $\OO$ are not empty.
Given an element $\LL\in\SSS$, the left regular representation $\rho$
induced by $\LL$, together with
the identity map, gives rise to a bijective 1-cocycle
$(\rho,\id)\in\OO$. Conversely,
for each $(f,q)\in \OO$, the multiplication
\begin{equation}\label{eq:3.2}
x\cdot y=q^{-1}(f(x)q(y)),\;\; \forall x,y\in \g,
\end{equation}
gives rise to an LSSA $(\LL,\cdot)\in\SSS$.
We denote the maps defined above by
$$\Psi: \SSS\ra \OO \quad {\text{and}} \quad   \Phi: \OO\ra \SSS.$$

To classify LSSAs  up to isomorphism, we introduce the notion of
quasi-equivalence between  1-cocycles so that $\Psi$ and $\Phi$ above
induce a bijection between isomorphic classes $\SSS/\cong$  in $\SSS$
and  quasi-equivalent classes $\OO/\simeq$ in $\OO$.

\begin{defn}{\rm Let $\g$ be a Lie superalgebra and  $f_i:\g\ra
\gl(V_i), i=1,2$, be two  representations of $\g$. Two  1-cocycles
$(f_1,q_1)$ and $(f_2,q_2)$ on  $\g$ are called \textit{equivalent},
denoted by $(f_1,q_1)\cong(f_2,q_2)$, if there exists a linear
isomorphism  $\varphi: V_2\ra V_1$ such that
$$f_2(x)=\varphi^{-1}\circ f_1(x)\circ \varphi \quad 
\textrm{ and } \quad q_2=\varphi^{-1}\circ q_1,\;\;\forall x\in\g.$$ We
say that two  1-cocycles $(f_1,q_1)$ and $(f_2,q_2)$ on $\g$ are
\textit{quasi-equivalent}, denoted by $(f_1,q_1)\simeq(f_2,q_2)$, 
if there exists an automorphism $T$ of  $\g$ such that $(f_1\circ T,
q_1\circ T)$ and $(f_2, q_2)$  are equivalent.
}\end{defn}

\begin{prop}\label{lem:3.3}
The maps $\Psi$ and $\Phi$ induce a bijection between the set
$\SSS/\cong$ of isomorphic classes of LSSAs on a Lie superalgebra $\g$
and the set  $\OO/\simeq$ of quasi-equivalent classes of bijective
1-cocycles on $\g$.
\end{prop}

\begin{proof}
Suppose that $T:\LL_2 \ra \LL_1$ is an isomorphism of  two LSSAs
$(\LL_1, \cdot_1)$ and $(\LL_2, \cdot_2)$   on  $\g$. Then $T$ is also
an automorphism  of the Lie superalgebra $\g$. %Set $\varphi=T$.  
 We have
$$\rho_{\LL_2}(x)(y)=x\cdot_2 y=T^{-1}(T(x)\cdot_1 T(y))=(T^{-1}\circ
(\rho_{\LL_1}\circ T)(x)\circ T)(y),\;\;\forall x,y\in\g.$$
Together with the  fact that $\id=T^{-1}\circ \id\circ T$, we deduce
that $(\rho_{\LL_1},\id)$ and $(\rho_{\LL_2},\id)$ are
quasi-equivalent. Hence isomorphic LSSAs are mapped to quasi-equivalent
bijective 1-cocycles.

Conversely, suppose that $(f_1,q_{1})$ and $(f_2,q_{2})$ are two
quasi-equivalent bijective 1-cocycles on $\g$. Then there exists a
linear isomorphism $\varphi$   and an automorphism $T$ of   $\g$ such
that $f_2(x)=\varphi^{-1}\circ(f_1\circ T)(x)\circ \varphi$ and
$q_{2}=\varphi^{-1}\circ q_{1}\circ T$ for all $x\in\g$.
Let  $(\LL_1,\cdot_1)$ and
$(\LL_2,\cdot_2)$ be the corresponding LSSAs induced from 
$(f_1,q_{1})$ and $(f_2,q_{2})$ by Eq. (\ref{eq:3.2}), respectively.
Then
$$T(x\cdot_2 y)=T(q_{2}^{-1}(f_2(x)q_2(y)))=q_{1}^{-1}
(f_1(T(x))q_1(T(y)))=T(x)\cdot_1  T(y)$$  
for all $x,y\in\g.$ Hence quasi-equivalent  bijective 1-cocycles  are
mapped to isomorphic LSSAs.

Note that $\Phi\circ \Psi(\overline{\LL})=\overline{\LL}$ for all
$\overline{\LL}\in\SSS/\cong$. For all $\overline{(f,q)}\in
\OO/\simeq$, we see that $\Psi\circ\Phi
(\overline{(f,q)})=\Psi(\overline{\LL})=\overline{(\rho_\LL,\id)}$.
Since $\rho_\LL(x)y=q^{-1}(f(x)q(y))$ for all $x,y\in \g$ and
$\id=q^{-1}\circ q$, we have  $(\rho_\LL,\id)$ is quasi-equivalent to
$(f,q)$. Therefore, $\Phi$ and $\Psi$ induce mutually inverse
bijections between $\SSS/\cong$ and  $\OO/\simeq$.
\end{proof}

\subsection{Evaluation maps} \label{sec2.3}
Given a representation $f:\g\ra\gl(V)$ of a Lie superalgebra $\g$ and
an element
$a\in V_{\0}$,  the map $\ev_a:\g\ra V$ defined by $\ev_a(x)=f(x)a$ for
all $x\in\g$ is called the \textit{evaluation map} of $\g$ associated
with  $f$ at the point $a$. It is immediate that $(f,\ev_a)$ is a
1-cocycle on $\g$. Evaluation maps are very useful in establishing
isomorphisms of LSSAs on $\g$; see Propositions  \ref{prop:3.6} and 
\ref{prop:3.7} below.

\begin{lem}\label{lem:3.5}
Let $f:\g\ra \gl(V)$ be a representation of a Lie superalgebra $\g$. If
there exist $a,b \in V_{\bar 0}$ satisfying $f(\g)a=f(\g)b=V$, then
$(f, \ev_a)$ and $(f, \ev_b)$ are quasi-equivalent.
\end{lem}

\begin{proof}
It is clear that $f(\g_{\bar 0})a= f(\g_{\bar 0})b= V_{\bar 0}$. Let 
$G_\0$ be the simply connected algebraic group with Lie algebra $\g_\0$
and let $F: G_{\0} \ra GL(V_{\0})$ be the representation of $G_\0$ with
$dF=f|_{\g_{\0}}$. Then both $F(G_\0)a$ and $F(G_\0)b$ are open in $V$
and hence $F(G_\0)a = F(G_\0)b$.
Choose $t\in G_{\bar 0}$ such that $F(t)a=b$ and
define $T\in \aut(\g)$ by $T=\Ad_t$ and $\varphi: V\ra V$ by
$\varphi=F(t)$.  Then $f\circ T:\g \ra\gl(V) $ sends every element
$x\in\g$ to $\varphi\circ f(x)\circ \varphi^{-1}$.
Hence, $f(x)=\varphi^{-1}\circ (f\circ T)(x)\circ \varphi$. Further,
$(\ev_b \circ T)(x)=\ev_b(T(x))=(f(T(x))(b)=(\varphi\circ f(x)\circ
\varphi^{-1})(b)=(\varphi\circ f(x))(a)= \varphi(\ev_a(x))$ for all
$x\in \g$ and hence $\ev_a=\varphi^{-1} \circ(\ev_b\circ T)$, as
desired.
\end{proof}

\begin{prop}\label{prop:3.6}
Let $f:\g\ra \gl(V)$ be a representation of a Lie superalgebra $\g$
with $\dim \g_{\gamma}=\dim V_{\gamma}$  for  $\gamma\in \Z_2$.
 If there exists an element  $a\in V_{\bar 0}$ such that the
 evaluation map $\ev_a$ is bijective and bijective 1-cocycles and
 bijective evaluation maps associated with $f$ coincide, then there
 exists  a unique  LSSA up to isomorphism on $\g$ associated with  $f$.
\end{prop}

\begin{proof}
We denote by $\LL$ the  LSSA on $\g$ given by $(f, \ev_a)$.
Suppose there exists another LSSA $\LL'$ on $\g$ given by the bijective
evaluation map $(f, \ev_b)$ with $b\in V_{\bar 0}$. Then
$\ev_a(\g)=\ev_b(\g)=V$, that is, $f(\g)a=f(\g)b=V$. Lemma
\ref{lem:3.5} implies that
$(f, \ev_a)$ and $(f, \ev_b)$ are quasi-equivalent. It follows from
Proposition \ref{lem:3.3} that $\LL$ is isomorphic to $\LL'$.
\end{proof}

\begin{prop}\label{prop:3.7}
Let $f_i:\g\ra \gl(V_i), i=1,2,$ be two quasi-equivalent
representations of $\g$, i.e.,
there exist $T\in\aut(\g)$ and an isomorphism $\varphi :V_2\ra V_1$
such that $f_2(x)=\varphi^{-1}\circ (f_1\circ T)(x)\circ \varphi$ 
for all $x\in\g$. Assume further that
$\dim(V_{i})_\gamma=\dim\g_\gamma$ for $i=1,2, \gamma\in \Z_2$ and that
bijective 1-cocycles and bijective evaluation maps associated with
$f_i$  coincide for each $i = 1,2$. Then we have
\begin{enumerate}
\item[(1)] if  there exists a bijective evaluation map associated with
one of them, then there exists  a bijective evaluation map   associated
with  the other one;  
\item[(2)] LSSAs associated with $f_1$ and $f_2$ are isomorphic.
\end{enumerate}
\end{prop}

\begin{proof}
(1) Suppose there exists an element $b \in (V_2)_{\bar 0}$  such that
the evaluation map $\ev_b$ associated with $f_2$ is bijective. Then
$\ev_b(\g)=f_2(\g)b=V_2$, and hence $V_2=(\varphi^{-1}\circ f_1(
T(\g))\circ \varphi)b=(\varphi^{-1}\circ f_1(\g)\circ \varphi)b$. Let
$a:=\varphi(b)\in (V_1)_{\bar 0}$. Then $\ev_a(\g)=f_1(\g)a=\varphi
(V_2)=V_1$, i.e., $\ev_a$ is surjective. Note that
$\dim V_1=\dim \g$. Then $\ev_a$ is injective and hence there exists
a bijective evaluation map $\ev_a$ associated with $f_1$.

(2) It follows from Proposition \ref{prop:3.6} that  there exists a
unique LSSA associated with each $f_i, i=1,2$.   We denote by $\LL$ and
$\LL'$ the  LSSAs on $\g$ given by $(f_1, \ev_a)$ and $(f_2, \ev_b)$,
respectively.
Since  $\ev_b(x)=f_2(x)b=(\varphi^{-1}\circ f_1(T(x)))a=(\varphi^{-1}
\circ \ev_a \circ T)(x)$, we see that $(f_1, \ev_a)$ and $(f_2, \ev_b)$
are quasi-equivalent. By Proposition \ref{lem:3.3}, we conclude that
$\LL$ and $\LL'$ are isomorphic.
\end{proof}

%%%%%%%%%%%%%%%%%%%%%%%%%%%%%%%%%%%%%%%%%%%SECTION 1
\section{Representations of $\sll(m|n)$} \label{sec3}
\setcounter{equation}{0}
\renewcommand{\theequation}
{3.\arabic{equation}}
\setcounter{theorem}{0}
\renewcommand{\thetheorem}
{3.\arabic{theorem}}

We present some preliminaries and calculations on the representations of
the Lie superalgebras  $\sll(m|n)$. For more details, see
\cite{Kac1977, Kac1978}.

\subsection{Definitions} \label{sec3.1}
Let  $\gl(m|n)$ be the space of $(m+n)\times (m+n)$ matrices. We write
an element $X \in \gl(m|n)$ in a block-diagonal form
$X = \left(
\begin{matrix}
X_1&X_2\\
X_3&X_4
\end{matrix}
\right)$, where $X_1, X_2, X_3, X_4$ are matrices of sizes  
$(m\times m)$, $(m\times n)$, $(n\times m)$, and $(n\times n)$
respectively.
Setting
\[\gl(m|n)_{\0} =\{ \left(
\begin{matrix}
X_1&0\\
0&X_4
\end{matrix}
\right)\} \quad \quad {\text{and}} \quad \quad \gl(m|n)_{\w} =\{ \left(
\begin{matrix}
0&X_2\\
X_3&0
\end{matrix}
\right)\} \]
endows $\gl(m|n)$ with a $\Z_2$-grading. The corresponding
supercommutator defined by
\[[X, Y] = XY - (-1)^{\gamma\,' \gamma\,''} YX\]
where $X \in \gl(m|n)_{\gamma\, '}$ and $Y \in \gl(m|n)_{\gamma\,''}$ turns
$\gl(m|n)$ into a Lie superalgebra.

Furthermore $\gl(m|n)$ admits a $\Z$-grading 
\[\gl(m|n) = \gl(m|n)_{-1} \oplus \gl(m|n)_{0} \oplus \gl(m|n)_{1}\]
defined by 
\[ \gl(m|n)_{-1} = \{\left(
\begin{matrix}
0&0\\
X_3&0
\end{matrix}
\right)\}, \quad \quad \gl(m|n)_{0}  = \{ \left(
\begin{matrix}
X_1&0\\
0&X_4
\end{matrix}
\right)\}, \quad \quad  \gl(m|n)_{1}  = \{\left(
\begin{matrix}
0&X_2\\
0&0
\end{matrix}
\right)\}.\]
The two gradings are compatible, i.e., 
\[\gl(m|n)_{\0} = \gl(m|n)_0 \quad \quad {\text{and}} \quad \quad
\gl(m|n)_{\w} = \gl(m|n)_{-1} \oplus \gl(m|n)_1.\]

The supertrace of $X \in \gl(m|n)$ is defined as
$\mathrm{str}(X)=\mathrm{tr}(X_1)-\mathrm{tr}(X_4)$.
 The special linear Lie superalgebra $\sll(m|n)$ is the subalgebra of
 $\gl(m|n)$ of traceless matrices:   
 \[\sll(m|n)=\{X\in\gl(m|n)\,|\,\mathrm{str}(X)=0\}.\]
Clearly, $\sll(m|n)_{\0}=\sll_m \oplus \sll_n\oplus \C$ is a reductive
Lie algebra.
Since $\sll(m|n)$ is isomorphic to $\sll(n|m)$, we always assume
that $m \geq  n\geq  1$. If $m\neq n$, then $\sll(m|n)$ is a simple Lie
superalgebra. On the other hand,
$\sll(m|m)$ has a one-dimensional centre $\C I_{2m}$ and the Lie
superalgebra $\psl(m|m) := \sll(m|m)/{\C I_{2m}}$ is simple; 
its even part $\psl(m|m)_{\0} \cong \sll_m \oplus \sll_m$ is
semisimple.

For the rest of the paper $\g$ will be  the Lie superalgebra
$\gl(m|n)$ or $\sll(m|n)$. 
If $\theta$ is an automorphism of $\g$ and $\rho: \g \to \gl(V)$ is a
representation of $\g$, the $\theta$-twist $V^\theta$ of $V$ 
is the module corresponding to the representation $\rho \circ \theta$.
If $\theta$ is an inner automorphism, then $V^\theta \cong V$.
The supertranspose of $X = \left(
\begin{matrix}
X_1&X_2\\
X_3&X_4
\end{matrix}
\right) \in \g$ is the matrix $X^{st} := \left(
\begin{matrix}
X_1^t&X_3^t\\
-X_2^t&X_4^t
\end{matrix}
\right)$. The supertranspose $st$ is an antiautomorphism of $\g$ while
$-st$ is an automorphism which is not an inner automorphism. 
If $V$ is an irreducible $\g$-module then $V^{-st} \cong V^*$. 
However, $V^{-st} \not \cong V^*$ in general.

  \subsection{Roots, positive roots} \label{sec3.2}
  Denote the subalgebra of $\g$ of diagonal matrices by $\h$. As usual,
  we denote by
  $E_{ij}$, $1 \leq i,j, \leq m+n$ the elementary matrix, i.e., the
  matrix with $1$ in position $(i,j)$ and zeroes elsewhere. 
  If $\g = \gl(m|n)$, let $\{\epsilon_1, \ldots, \epsilon_m, \delta_1,
  \ldots, \delta_n\}$ be the basis of 
  $\h^*$ dual to the basis $\{E_{11}, \ldots, E_{mm}, E_{m+1,m+1},
  \ldots, E_{nn}\}$ of $\h$. 
 For $\g = \sll(m|n)$, slightly abusing notation, we denote the
 restriction of $\epsilon_i$ and $\delta_j$ to $\h$ by $\epsilon_i$ and
 $\delta_j$ as well. 
 Note that $\{\epsilon_1, \ldots, \epsilon_m, \delta_1, \ldots,
 \delta_n\}$ span $\h^*$ and satisfy the relation
 \[\epsilon_1 + \ldots + \epsilon_m = \delta_1 + \ldots + \delta_n.\]

The Lie superalgebra $\g$ admits a root decomposition 
\[\g = \h \oplus \left( \oplus_{\alpha \in \Delta} \,
\g^\alpha\right),\]
where, for any $\alpha \in \h^*$, 
\[\g^\alpha = \{X \in \g \, | \, [h, X] = \alpha(h) X 
{\text{ for every }} h \in \h\} \quad 
{\text{and}} \quad
\Delta = \{\alpha \in \h^* \backslash \{0\} \, |\, 
\g^\alpha \neq 0\}.\]
The elements of $\Delta$ are called roots of $\g$. Furthermore,
the decomposition $\g = \g_{\0} \oplus \g_{\w}$ induces the decomposition
of $\Delta$  as
\[\Delta = \Delta_\0 \cup \Delta_{\w}, \quad {\text{where}} \quad
\Delta_\gamma = 
\{\alpha \in \Delta \, | \, \g^\alpha \subset \g_\gamma\}, \gamma\in \Z_2.\]
Explicitly,
\[\Delta_{\0} = \{ \epsilon_i - \epsilon_j \, | \, 1 \leq i \neq j 
\leq m\} \cup \{\delta_i - \delta_j \, | \, 1 \leq i \neq j \leq n\},\]
and
\[\Delta_{\w} = \{ \pm(\epsilon_i  - \delta_j) \, | \, 1 \leq i  
\leq m, 1 \leq j \leq n\}.\] 
All root spaces are one-dimensional and spanned by elementary matrices.
Namely, $\g^{\epsilon_i - \epsilon_j}$, $\g^{\delta_i - \delta_j}$, 
$\g^{\epsilon_i - \delta_j}$, $\g^{\delta_i - \epsilon_j}$ are spanned
by $E_{ij}, E_{m+i, m+j}, E_{i, m+j}, E_{m+i, j}$ respectively.

Let $\Delta_{-1}$ and $\Delta_1$ denote the roots of $\g_{-1}$ and
$\g_1$ respectively and let 
\[\Delta_{\0} = \Delta_{\0}^+ \cup \Delta_{\0}^-\]
be the triangular decomposition of $\Delta_{\0}$ defined by
$\Delta_{\0}^\pm = \{\pm(\epsilon_i - \epsilon_j) \, | \, 
1 \leq i < j \leq m\}$. 
Fixing the set $\Delta^+ = \Delta_{\0}^+ \cup \Delta_1$ of positive
roots of $\g$, we denote the corresponding Borel subalgebra of $\g$ 
by $\gb$:
\[\gb =  \h \oplus \left( \oplus_{\alpha \in \Delta^+} \,
\g^\alpha\right).\]
Note that $\gb_{\0} = \gb \cap \g_{\0}$ is a Borel subalgebra of
$\g_{\0}$ with roots $\Delta_{\0}^+$. The $\Z$-grading of $\g$ defines
the parabolic 
subalgebra $\p = \g_0 \oplus \g_1$ with roots $\Delta_{\0} \cup
\Delta_{\w}^+$.

\subsection{Representations} \label{sec3.3}
Given a (finite-dimensional) $\g_{\0}$-module $L$, setting $\g_1 \cdot
L = 0$, we turn it into a $\p$-module and define the corresponding 
parabolically induced module $K(L)$ by
\[ K(L):=\mathrm{Ind}^{\g}_{\p} \, L \,\,\cong \wedge^{\cdot}
(\g_{-1})\otimes_{\C}L.\]
Here $\wedge^{\cdot} (\g_{-1})$ denotes the exterior algebra of the
vector space $\g_{-1}$ and
the isomorphism is an isomorphism of $\g_{-1}$-modules.

 We define a symmetric bilinear form on $\h^*$ by $(\epsilon_i,
 \epsilon_i)= 1, (\delta_j, \delta_j)=-1$ and setting all other
 pairings between 
 elements $\epsilon_1, \ldots, \epsilon_m, \delta_1, \ldots, \delta_n$
 to be equal to zero.
 A   weight $\lambda\in \h^*$ is said to be integral if 
 $(\lambda, \beta)\in \Z$ for all roots $\beta \in \Delta_{\0}$, 
 and dominant if $2\frac{(\lambda, \beta)}{(\beta, \beta)}\geq 0$ 
 for all $\beta \in \Delta^+_{\0}$. 
 We denote by $X^+$ the set of dominant integral weights in $\h^*$.
It parametrizes the isomorphism classes of irreducible
finite-dimensional $\g_{\0}$-modules. For a given $\lambda \in X^+$,
the corresponding $\g_{\0}$-module is denoted by $L(\lambda)$. The Kac
module $K(\lambda)$ is simply  $K(L(\lambda))$. It admits a unique
proper maximal submodule $I(\lambda)$ and, respectively, a unique irreducible
quotient $V(\lambda) = K(\lambda)/I(\lambda)$. Every irreducible 
$\g$-module is isomorphic to
$V(\lambda)$ or $\Pi V(\lambda)$
for a unique $\lambda \in X^+$, where $\Pi$ is the parity change
functor.  Note that, whenever using the notations $K(\lambda)$ and 
$V(\lambda)$, we assume that the highest weight space is even.

%Whenever we use the notation $\Pi V(\lambda)$, we will assume
%that the highest weight vector of $V(\lambda)$ is even (and hence the 
%highest weight vector of $\Pi V(\lambda)$ is odd). However, 
%sometimes we will use $V(\lambda)$ to denote the highest weight module 
%with an odd highest weight vector. In all such cases, we will indicate 
%the parity of the highest weight vector. 

The weight $\lambda\in X^+$  is called typical if
$(\lambda+\rho, \alpha)\neq0$ for all $\alpha\in\Delta^+_{\w}$, where 
\[\rho = \rho_{\0} - \rho_{\w},\]
with
\[\rho_{\0} = \frac{1}{2} \sum_{\alpha \in \Delta_{\0}^+} \, 
\alpha \quad \quad {\text{and}} \quad \quad 
\rho_{\w} = \frac{1}{2} \sum_{\alpha \in \Delta_{\w}^+} \, \alpha.\]
If $\lambda$ is not typical, it is called atypical. The modules
$K(\lambda)$ and $V(\lambda)$ are called typical (respectively,
atypical) if the corresponding weight $\lambda$ is typical (respectively,
atypical). Note that $K(\lambda)$ is irreducible if and
only if it is typical; otherwise $K(\lambda)$ is indecomposable but
reducible.

The degree of atypicality of $\lambda$ (and of the respective modules
$K(\lambda)$ and $V(\lambda)$) is defined as  the number of distinct
elements $\alpha\in\Delta^+_{\w}$ for which $(\lambda+\rho, \alpha)=0$.
If there exists one and only one such $\alpha\in\Delta^+_{\w}$,
$\lambda$  is called singly atypical.
Note that any dominant integral  weight of $\sll(m|1), m\geq2$ is
either typical or singly atypical (\cite[Lemma 3.2.1]{Ger1998}).
Germony studied singly atypical representations in \cite{Ger1998} and
we will rely on the results therein.

If $\lambda$ is a singly atypical weight, then $K(\lambda)$ contains a
unique proper submodule which is irreducible. Let $T^- \lambda$ denote 
the highest weight of the unique proper submodule  of
$K(\lambda)$. In other words, $K(\lambda)$ is a non-split extension of
$V(\lambda)$ by $V(T^- \lambda)$ or by $\Pi V(T^- \lambda)$, 
depending on the parity of the highest weight space of the proper submodule
of $K(\lambda)$. The operator $T^-$ has an inverse
denoted by $T^+$.  The extensions among simple $\g$-modules is
described in \cite[Proposition 6.1.2]{Ger1998}:

\begin{prop}\label{prop2.7}
Let $\lambda, \mu \in X^+$  be dominant integral  weights.
\begin{enumerate}
  \item If $\lambda$ is typical, then
    \begin{equation}
  \mathrm{Ext}^1(V(\lambda), V(\mu))=
\begin{cases}
\C& \text{if  $\lambda=\mu$},\\
0& \text{otherwise}.
\end{cases}
\end{equation}
  \item If $\lambda$ is singly atypical, then
   \begin{equation}
\dim  \mathrm{Ext}^1(V(\lambda), V(\mu))=
\begin{cases}
1& \text{if  $\mu \in \{T^+\lambda, T^-\lambda\}$},\\
0& \text{otherwise}.
\end{cases}
\end{equation}
\end{enumerate}
\end{prop}

As mentioned above, $K(\lambda)$ is a non-split extension of
$V(\lambda)$ by $V(T^- \lambda)$ or by $\Pi V(T^- \lambda)$. 
To describe a non-split extension of $V(\lambda)$ 
by $V(T^+ \lambda)$ or by 
$\Pi V(T^+ \lambda)$, we introduce the following notation. Given 
$\mu \in X^+$, the opposite Kac module $K'(\mu)$ is the module
\[ K'(\mu):=\mathrm{Ind}^{\g}_{\g_{-1} \oplus \g_0} \, L(\mu)\]
and the weight $\mu'$ is defined as the unique element of $X^+$ for
which $V(\mu)$ is a quotient of $K'(\mu')$.
The module $K'(\lambda')$ is a non-split extension of $V(\lambda)$ by
$V(T^+ \lambda)$ or by $\Pi V(T^+ \lambda)$. 
Thus $K(\lambda)$ and $K'(\lambda')$ provide examples
for Proposition \ref{prop2.7} (2).

Next we provide an example for 
Proposition \ref{prop2.7} (1). Let $\g = \sll(m|n)$ and denote 
by $\C^{(2)}$ the two-dimensional
$\g_\0$-module on which every element acts trivially, except
that a fixed nonzero central element $z$ acts via a nilpotent matrix of
order 2 and set 
\[K(\lambda)^{(2)}= K(L(\lambda)\otimes \C^{(2)}) =
\Ind^\g_{\p} (L(\lambda)\otimes \C^{(2)}).\]
It is clear that $K(\lambda)^{(2)}$ is a non-split extension of
$K(\lambda)$ by itself. If $\lambda$ is singly atypical, the structure
of $K(\lambda)^{(2)}$ is described in \cite[Lemma 6.1.1]{Ger1998}:

\begin{lem}\label{lem2.1}
Let $\lambda$ be a singly atypical dominant integral  weight. The
module $K(\lambda)^{(2)}$ is uniserial with composition factors
(listed from top to socle) $V(\lambda), V(T^-\lambda), V(\lambda)$, and
$V(T^-\lambda)$ or 
$V(\lambda), \Pi V(T^-\lambda), V(\lambda)$, and $\Pi V(T^-\lambda)$.
\end{lem}

\subsection{Irreducible representations of $\sll(2|1)$} \label{sec3.4}
Let $\g = \sll(2|1)$. Fix the elements $h = E_{11} - E_{22}$ and $z 
=E_{11} + E_{22} + 2E_{33}$ which form a basis of the Cartan subalgebra
$\h$ of $\g$.
Let $V(i, k)$ for $(i,k) \in \Z_{\geq 0} \times \C$ denote  the
irreducible $\g$-module with highest weight $\lambda$ defined by 
$\lambda(h) = i$ and $\lambda(z) = k$.
The corresponding Kac module and irreducuble  module are denoted 
respectively by $K(i,k)$ and $V(i,k)$.  Denote
by $S_i$ the irreducible $\sll_2$-module of dimension $i+1$; by
convention, 
$S_{-1} = 0$.
The following proposition describes the modules $V(i,k)$ and the
extensions among them.

\begin{prop}\label{prop3.10} Let $(i, k) \in \Z_{\geq 0} \times \C$.

\begin{enumerate}
\item $(i,k)$ is typical if and only if $k \not \in \{i, -i-2\}$. 
\item If $k \not \in \{i, -i-2\}$, $V(i,k) = K(i,k)$ is of dimension
$2(i+1)|2(i+1)$ and, as an $\sll_2$-module, 
\[V(i,k)_{\0} \cong S_i \oplus S_i, \quad \quad \quad V(i,k)_{\w} 
\cong S_{i-1} \oplus S_{i+1}.\]
Moreover, $z$ acts on $V(i,k)_{\w}$ as multiplication by $k+1$ and
the decomposition above can be chosen so that $z$ acts on one 
of the copies of $V(i,k)_{\0}$ as multiplication by $k$ and on
the other one as multiplication by $k+2$. 
\item The module $V(i,i)$ is of dimension $i+1|i$ and, as an
$\sll_2$-module, 
\[V(i,i)_{\0} \cong S_i, \quad \quad \quad V(i,i)_{\w} \cong S_{i-1}.\]
Moreover, $z$ acts on $V(i,i)_{\0}$ as multiplication by $k$ and on
$V(i,k)_{\w}$ as multiplication by $k+1$.

\noindent
The module $V(i,-i-2)$ is of dimension $i+1|i+2$ and, as an
$\sll_2$-module, 
\[V(i,i)_{\0} \cong S_i, \quad \quad \quad V(i,i)_{\w} \cong S_{i+1}.\]
Moreover, $z$ acts on $V(i,i)_{\0}$ as multiplication by $k$ and on $V(i,k)_{\w}$ as
multiplication by $k+1$.
\item The operator $T^-$ acts on atypical weights as 
\[T^-(i,k) = \left\{ \begin{array}{ccccc}
 (i+1, i+1) & \quad & {\text{if}} & \quad& k = i\\
 (i-1, -i-1) & \quad& {\text{if}} & \quad& i > 0,\, k = -i-2\\
(0, 0)& \quad & {\text{if}} & \quad&  i=0,\, k = -2. \end{array}
\right.\]
\end{enumerate}
\end{prop}

\begin{proof} Since $\rho = - \epsilon_2 + \delta$, for 
$\lambda =(i,k)$, we have
\[(\lambda + \rho, \epsilon_2 - \delta) = \frac{k-i}{2} \quad 
\quad {\text{and}} \quad \quad (\lambda + \rho, \epsilon_1 - \delta) 
=  \frac{k+i}{2} +1,\]
proving (1). The remaining statements follow from an easy calculation
using the explicit $\gl_2$-structure  of the module $K(i,k)$. We leave
these to the reader.
\end{proof}

It will be convenient to index the atypical weights of $\sll(2|1)$ by $\Z$.
Namely, set 
\[\lambda_i := \left\{ \begin{array}{ccccc} (i,i) & \quad & {\text{if}} &
\quad & i \geq 0\\
(-i-1,i-1) & \quad & {\text{if}} & \quad & i < 0. \end{array}
\right.\]
In this notation $T^-(\lambda_i) = \lambda_{i+1}$ for any $i \in \Z$. Setting
$V_i := V(\lambda_i)$, we conclude that 
there are non-trivial extensions between $V_i$ and $V_j$
if and only if $|i-j| = 1$. Taking into account parity, the non-trivial
extensions between 
$V_{i-1}$ and $V_i$ require that these modules are taken in different
parities except for $i=0$ when the parities have to be the same.

%In what follows, we will indicate which parity is used by choosing the
%dimension of one of these modules, e.g., writing $\dim V_3 = 4|3$ means
%that the highest weight vector is chosen to be even and writing 
%$\dim V_3 = 3|4$
%means that the highest weight vector is chosen to be odd. 
%In the latter case, we
%need to switch the respective $V_{\0}$ and $V_{\w}$ in Proposition
%\ref{prop3.10}.

We complete the discussion of irreducible $\sll(2|1)$-modules by
describing their twists by the automorphism $-st$.

\begin{prop} \label{prop3.11}
Let $V$ be an irreducible $\sll(2|1)$-module. Then
\[V^{-st} \cong \left\{ \begin{array}{cccccc} K(i, -k -2) & \quad &
{\text{if}} & \quad & V = K(i,k)&{\text{is typical}}\\
\Pi V_{-i}& \quad & {\text{if}} & \quad & V = V_i, \ i \neq 0 & \\ %{\text{is atypical}}.
V_0 & \quad& {\text{if}} & \quad & V = V_0\, . & 
\end{array} \right.
\]
\end{prop}

\begin{proof} 
If $w$ is the lowest weight vector of $V$ and the weight of $w$ is
$\mu$, then $w$ is the highest weight vector of $V^{-st}$ of weight
$-\mu$. An explicit calculation which determines the lowest weighs of
the irreducible representations of $\sll(2|1)$ completes the proof. We
omit this calculation here.
\end{proof}

\subsection{Small irreducible modules of $\sll(m|1)$ for $m \geq 3$}
\label{sec3.5}

Let $\C^m$ be the standard $m$-dimensional module of
$\sll_m$ and $\C^{m|n}$  the standard $m|n$-dimensional module $\sll(m|n)$.
As $\sll_m$-modules, we have $\sll(m|1)_{-1}\cong (\C^m)^*$ and
$\sll(m|1)_{1}\cong \C^m$, where $(\C^m)^*$ is the dual module of the
standard module $\C^{m}$ of $\sll_m$.
In this section, we denote by $\tr$ the $1$-dimensional trivial module
of $\sll_m$ and by $\ttr$ -- the $0|1$-dimensional trivial module of
$\sll(m|1)$.

\begin{prop}\label{lem:4.4}
 Let $V$ be an irreducible $\sll(m|1)$-module with $m\geq 3$.
 \begin{enumerate}
   \item If $V$ is purely odd, i.e., $V_{\0}=0$ (respectively, purely
   even), then $V$ is the $0|1$-dimensional (respectively,
   $1|0$-dimensional) trivial module.
   \item Assume that $1 < \dim V \leq m^2|2m$, i.e., $V$ is non-trivial and
   $\dim V_{\0} \leq m^2$ and 
   $\dim V_{\w} \leq 2m$. Furthermore, assume that, as an $\sll_m$-module, 
   $V_{\0}$ is isomorphic to a direct sum of copies of $\C^m$. Then
   $V$ is isomorphic to one of the following modules:
        \begin{enumerate}
          \item $\C^{m|1}$ {\rm ;}
          \item $\Pi \wedge^2(\C^{m|1})$ for $m = 3, 4$ {\rm ;}
          \item $\Pi S^2(\C^{3|1})$ for $m = 3$.
        \end{enumerate}
 \end{enumerate}
\end{prop}

\begin{proof}
(1) If $V$ is purely even or purely odd, then $\sll(m|1)_{\w} 
\cdot V = 0$. Since $\sll(m|1)_{\w}$ generates $\sll(m|1)$, we 
conclude that $\sll(m|1) \cdot V = 0$. The irreducibility of $V$
implies $\dim V=0|1$ or  $\dim V=1|0$.

(2) We  assume that $V$ is induced from an irreducible $\gl_m$-module
$L(\lambda)$ with the highest weight $\lambda$, then, up to parity,
$V$ is isomorphic to the Kac module $K(\lambda)$ or a quotient
$K(\lambda)/I(\lambda)$ of
$K(\lambda)$. Since $L(\lambda)$ is an irreducible $\gl_m$-module, it
is irreducible as an $\sll_m$-module.  Our proof will be separated
into two cases: (I) $L(\lambda)\subseteq V_{\0}$  and (II)
$L(\lambda)\subseteq V_{\w}$.

Case (I).
Since, as an $\sll_m$-module, $V_{\0}$ is isomorphic to a direct sum 
of copies of $\C^m$, we have $L(\lambda)\cong \C^m$. 
Thus $\lambda = \epsilon_1 + \mu \delta$,
where $\delta = \delta_1$, cf. Section \ref{sec3.2}.
 We  analyze the module structure of Kac module 
 $K(\lambda)=\wedge^{\cdot}(\sll(m|1)_{-1})\otimes \C^m$. The $\Z$-grading
 on $\wedge^{\cdot}(\sll(m|1)_{-1})$ induces a grading
 \[K(\lambda) = \oplus_{i = 0}^m \, K_i\, .\]
The component $K_1$ is
contained in the odd part of $K(\lambda)$ and, as an $\sll_m$-module, it 
is isomorphic to the
direct sum of the trivial module and  the adjoint module of dimension $m^2-1$.
Since $m^2-1 \geq 2m$ for $m\geq 3$,  the adjoint $\sll_m$-module must be
contained in the unique submodule $I(\lambda)$
of $K(\lambda)$. Consequently, $\mu = 0$ and $\lambda = \epsilon_1$,
proving that $V(\lambda) \cong \C^{m|1}$ with $V_{\0}\cong \C^m$ and 
$V_{\w}\cong \tr$
as $\sll_m$-modules.

 Case (II). Since $L(\lambda) \subset V_{\w}$, the inequalities 
 \[\dim L(\lambda) \leq \dim V_{\w} \leq 2m\]
 imply that  $L(\lambda)$, as an
 $\sll_m$-module, is isomorphic to one of the following:
 \begin{enumerate}
     \item[] $\tr,\C^m,(\C^m)^* \quad$ for $m \geq 6$;
     \item[] $\tr,\C^m,(\C^m)^*, \wedge^2(\C^m), \wedge^2((\C^m)^*) \quad$ 
     for $m = 4, 5$;
     \item[] $\tr,\C^3,(\C^3)^*, S^2(\C^3), S^2((\C^3)^*) \quad$ 
     for $m = 3$.
 \end{enumerate}
 Here $\wedge^2(W)$ and $S^2(W)$ denote respectively the second exterior and 
 symmetric powers of $W$.
 Below  we  consider each of these cases for $L(\lambda)$.

(i)  If $L(\lambda) \cong \tr$, then $\lambda = \mu \delta$.
Then, as an $\sll_m$-module, $K_1 \cong (\C^m)^*$. However $K_1$ is 
contained in the even part of $K(\lambda)$ which must be a sum of copies of 
$\C^m$. Since $\C^m \not \cong (\C^m)^*$, we conclude that 
$K_1 \subset I(\lambda)$ which leads to $\lambda = 0$. This contradicts the
assumption that $V$ is non-trivial.

(ii) If $L(\lambda) \cong \C^m$,
then $\lambda= \epsilon_1 + \mu \delta$. 
As in Case (I), $K_1$ is the direct sum of the
trivial  and the adjoint modules of $\sll_m$. These two modules are
not isomorphic to $\C^m$. Thus $K_1 \subset I(\lambda)$ and 
$V(\lambda) = K(\lambda)/I(\lambda) = K_0$. In particular, $V$ is purely
odd and, by (1), $V$ is trivial.

(iii) If $L(\lambda) \cong (\C^m)^*$,
then $\lambda= -\epsilon_m + \mu \delta$. 
Then, as an $\sll_m$-module, $K_1$ is the direct sum of 
$S^2((\C^m)^*)$ and $\wedge^2((\C^m)^*)$, neither of which is 
isomorphic to $\C^m$ if $m \geq 4$. Arguing as in (ii) above, we conclude
that, for $m \geq 4$, $V$ is purely odd and thus trivial.
When $m=3$, we observe that $\wedge^2((\C^3)^*)\cong \C^3$ but 
$S^2((\C^m)^*) \not \cong \C^3$.
Thus  $S^2((\C^m)^*) \subset I(\lambda)$, 
$\lambda = - \epsilon_3 + \delta = \epsilon_1 + \epsilon_2$,
and $V \cong \Pi \wedge^2(\C^{3|1})$. 

(iv) If $m = 5$ and $L(\lambda) \cong \wedge^2(\C^5)$, then
$\lambda = \epsilon_1 + \epsilon_2 + \mu \delta$. 
Then, as an $\sll_5$-module, $K_1$ is the direct sum of 
$\C^5$ and the $45$-dimensional module with highest weight 
$\epsilon_1+ \epsilon_2 - \epsilon_5$. The assumption $K_1 \subset I(\lambda)$
leads to a contradiction as $V$ would be purely odd. Alternatively, 
$K_1 \cap I(\lambda)$ equals the $45$-dimensional module above. 
Then $\lambda = \epsilon_1 + \epsilon_2$ and $V = \Pi \wedge^2(\C^{5|1})$.
However, $\dim (\Pi \wedge^2(\C^{5|1}))_{\w} = 11 > 10$, contradicting
the assumption on $V$.

(v) If $m = 5$ and $L(\lambda) \cong \wedge^2((\C^5)^*)$, then
$\lambda = -\epsilon_4 - \epsilon_5 + \mu \delta$. 
Then, as an $\sll_5$-module, $K_1$ is the direct sum of 
$\wedge^2(\C^5)$ and the $40$-dimensional module with highest weight 
$-\epsilon_4 - 2\epsilon_5$. The assumption on $V$ implies that 
$K_1 \subset I(\lambda)$ which 
leads to a contradiction as $V$ would be purely odd.

(vi) If $m = 4$ and $L(\lambda) \cong \wedge^2(\C^4) \cong \wedge^2((\C^4)^*)$,
then $\lambda = \epsilon_1 + \epsilon_2 + \mu \delta$.
Then, as an $\sll_4$-module, $K_1$ is the direct sum of 
$\C^4$ and the $20$-dimensional module with highest weight 
$\epsilon_1 + \epsilon_2 - \epsilon_4$. As in (iv) above, we conclude that
$K_1 \cap I(\lambda)$ equals the $20$-dimensional module above. 
Then $\lambda = \epsilon_1 + \epsilon_2$ and 
$V = \Pi \wedge^2(\C^{4|1})$.

(vii) If $m = 3$ and $L(\lambda) \cong S^2(\C^3)$,
then $\lambda = 2\epsilon_1  + \mu \delta$.
Then, as an $\sll_3$-module, $K_1$ is the direct sum of 
$\C^3$ and the $15$-dimensional module with highest weight 
$2\epsilon_1 - \epsilon_3$. As in (vi) above, we conclude that
$V = \Pi S^2(\C^{3|1})$.

(viii) If $m = 3$ and $L(\lambda) \cong S^2((\C^3)^*)$,
then $\lambda = -2\epsilon_3  + \mu \delta$. Then, 
as an $\sll_3$-module, $K_1$ is the direct sum of 
the adjoint module and $S^3((\C^3)^*)$. Arguing as in (v) above, we reach
a contradiction.
\end{proof}

\begin{rem} Let $V$ be an irreducible $\sll(m|1)$-module for $m\geq 3$.
Assume that the dimension of $V$ is less than or equal to $m^2|2m$ and
$V_{\0}$, as an $\sll_m$-module, is isomorphic to a direct sum of 
copies of $(\C^m)^*$. Then $V^*$ satisfies the assumptions of 
Proposition \ref{lem:4.4} and hence $V$ is isomorphic to a module dual to 
one of the modules listed in Proposition \ref{lem:4.4}.
\end{rem}

We complete the discussion of $\sll(m|1)$-modules by recording some information
about the modules $\wedge^2(\C^{m|1})$ and $S^2(\C^{m|1})$. The proof is trivial
and we omit it here.

\begin{prop} \label{prop3.99} We have
\begin{enumerate}
    \item $\dim \wedge^2(\C^{m|1}) = \frac{m^2-m+2}{2}|m$ and 
    $\dim S^2(\C^{m|1}) = \frac{m^2+m}{2}|m$ {\rm ;}
    \item $(\wedge^2(\C^{m|1}))_{\0} \cong \wedge^2(\C^m) \oplus \tr$, 
    $(\wedge^2(\C^{m|1}))_{\w} \cong \C^m$ \\ and \\
    $(S^2(\C^{m|1}))_{\0} \cong S^2(\C^m)$, 
    $(S^2(\C^{m|1}))_{\w} \cong \C^m$
    as $\sll_m$-modules. 
\end{enumerate}
\end{prop}

\section{LSSAs on $\sll(m|1)$} \label{sec4}
\setcounter{equation}{0}
\renewcommand{\theequation}
{4.\arabic{equation}}
\setcounter{theorem}{0}
\renewcommand{\thetheorem}
{4.\arabic{theorem}}

The purpose of this section is to give a proof of Theorem \ref{thm1},
that is, we want to prove that there are no LSSAs on $\sll(m|1)$ for
$m\geq 3.$  By Lemma \ref{lem:3.3}, it suffices to show that the set
$\OO$ of bijective 1-cocycles of $\sll(m|1)$ is empty.

Throughout this section we assume $m \geq 3$. 
Let $P_{m}:=m\C^{m|1}\oplus m\, \widetilde{\tr}$ be the $\sll(m|1)$-module which is
the direct sum of $m$ copies of $\C^{m|1}$ and $m$ copies of
$\widetilde{\tr}$ and let $P_{m}^*$ be module dual to $P_{m}$.

\begin{prop}\label{prop:4.4}
Let $W$ be an $\sll(m|1)$-module of dimension $m^2|2m$ such that,  
as an $\sll_m$-module, $W_{\0}$ is isomorphic to the direct sum of
$m$ copies of $\C^m$ or $m$ copies of $(\C^m)^*$.

\begin{enumerate}
    \item If $m\geq4$, then $W$ is isomorphic to either $P_{m}$ or $P_{m}^*$.
    \item If $m=3$, then $W$ is isomorphic to one of $P_3,P_3^*,Q_3$, or $Q_3^*$,
    where \\ $Q_3= 2\C^{3|1}\oplus \Pi \wedge^2(\C^{3|1})$, see Proposition \ref{lem:4.4}, 
    and $Q_3^*$ is the module dual to $Q_3$.
\end{enumerate}
\end{prop}

\begin{proof} We will prove the proposition in the case when  
$W_{\0}$ is isomorphic to the direct sum of $m$ copies of $\C^m$. The case when
$W_{\0}$ is isomorphic to the direct sum of $m$ copies of $(\C^m)^*$ then
follows by duality.

If $0=W^0\subset W^1\subset W^2\subset\cdots\subset W^k=W$ is
a composition series of $W$, then
each $V^i:=W^i/W^{i-1}$ is an irreducible  $\sll(m|1)$-module for
$1\leq i\leq k$.
If $\dim V^i=a_i|b_i$, then $\dim W=m^2|m$ implies that $\sum_{i=1}^k a_i=m^2$
and $\sum_{i=1}^k b_i=2m$. Hence $\dim V^i\leq m^2|2m$ for $1 \leq i \leq k$ and
$V^i$ is isomorphic to $\ttr$ or one of the modules from Proposition \ref{lem:4.4}.
Combining Propositions \ref{lem:4.4} and \ref{prop3.99}, we also get 
$\dim(V^i)_{\0} \leq m$, $\dim (V^i)_{\w} \geq 1$, and $k \geq m$.

As an  $\sll_m$-module, $W_{\0}\cong \oplus_{i=1}^k (V^i)_{\0}$, implying that each
$(V^i)_{\0}$ itself is isomorphic to a (possibly empty) direct sum  of
copies of  $\C^m$  for $1\leq i\leq k$.
Proposition \ref{lem:4.4} implies that each $V^i$ is isomorphic to one of the 
following modules

\begin{enumerate}
     \item[] $\C^{m|1}$ or $\ttr \quad$ for $m \geq 5$;
     \item[] $\C^{4|1}$, $\ttr$, or $\Pi \wedge^2(\C^{4|1}) \quad$ for $m = 4$;
     \item[] $\C^{3|1}$, $\ttr$, $\Pi \wedge^2(\C^{3|1})$, or
     $\Pi S^2(\C^{3|1})\quad$ for $m = 3$.
 \end{enumerate}

First we note that, for $m = 4$, $V^i$ cannot be isomorphic to 
$\Pi \wedge^2(\C^{4|1})$. If, to the contrary, $V^i \cong \Pi \wedge^2(\C^{4|1})$
for some $i$, then
\[\dim W_{\w} = \dim (V^i)_{\w} + \sum_{j \neq i} \dim (V^j)_{\w} \geq 
7 + (k-1) \geq 7 + 3 > 8, \]
which contradicts the assumption on $W$. A similar argument shows that, 
for $m = 3$, $V^i$ cannot be isomorphic to $\Pi S^2(\C^{3|1})$. 
This proves that each $V^i$ 
is isomorphic to one of the modules $\C^{m|1}, \ttr$, or $\Pi \wedge^2(\C^{3|1})$.

Counting dimensions we conclude that either
\begin{enumerate}
    \item[(i)] $k = 2m$ and 
    $V^i \cong \C^{m|1}$ for $m$ values of $i$ and $V^i \cong \ttr$ for
    the remaining $m$ values of $i$
    
    or
    
    \item[(ii)] $m = k = 3$ and $V^i \cong \Pi \wedge^2(\C^{3|1})$ for one value of
    $i$ and $V^i \cong \C^{3|1}$ for two values of $i$.
\end{enumerate}
However, by Proposition \ref{prop2.7}, there are no non-trivial extensions
between $\C^{m|1}$ and $\ttr$ and there are no non-trivial 
extensions between 
$\Pi\wedge^2(\C^{3|1})$ and $\C^{3|1}$. Thus $W$ is isomorphic to $P_m$ or $Q_3$.
\end{proof}

\begin{lem}\label{lemma:4.1}
All bijective 1-cocycles of $\gl_m$ are bijective evaluation maps for
$m>1$.
\end{lem}

\begin{proof}
Let $q: \gl_m \longrightarrow V$ be a bijective  1-cocycle of $\gl_m$
associated with the representation $f:\gl_m \ra \gl(V)$.
Then $q$ induces a left-symmetric algebra structure on $\gl_m$, see 
\cite[Theorem 2.1]{Bai2009}. The results of Bauers, \cite{Bau1999},
imply that $q$ is a bijective evaluation map. Namely,
it follows from \cite[Section 2.2]{Bau1999} that  there exists an
\'{e}tale affine representations of $\gl_m$ with base point $0\in V$ and
evaluation map $ev_0=q$. Furthermore, 
all \'{e}tale affine representations of $\gl_m$ are linear,
see \cite[Propositions 5.1 and 2.2]{Bau1999} and hence all
bijective 1-cocycles of $\gl_m$ are bijective evaluation maps.
\end{proof}

\begin{lem}\label{lemma:4.2}
If $q$ is a bijective 1-cocycle of $\sll(m|1)$ associated with one of
the modules $P_{m}, P_{m}^*$, for $m \geq 3$, $Q_3$, or $Q_3^*$,
then $q$ is a bijective evaluation map.
\end{lem}

\begin{proof}
Let $q$ be a bijective 1-cocycle of $\sll(m|1)$ associated with the module $V$.
Since the even part of $\sll(m|1)$ is $\gl_m$, we
see that the restriction $q|_{\gl_m}$  is a bijective 1-cocycle of
$\gl_m$. Then $q|_{\gl_m}$ is a bijective evaluation map by Lemma
\ref{lemma:4.1}. Let $q|_{\gl_m}=\ev_a:\gl_m \ra V_{\0}$
for some point $a\in V_{\bar 0}$ and extend $\ev_a$ to the evaluation map
$\widetilde{\ev}_a:\sll(m|1)\ra V$ at the same point $a$ by setting
$\widetilde{\ev}_a(y)=y\cdot a$ for all $y\in \sll(m|1)$.  
Define $p:= q-\widetilde{\ev}_a$.
Then $p$ is a 1-cocycle of $\sll(m|1)$ associated with $V$ and $p|_{\gl_m}=0$.
It follows from Lemma \ref{lem:3.1} that   $p|_{\sll(m|1)_{\bar 1}}:
\sll(m|1)_{\bar 1}\ra V_{\w}$ is a homomorphism of $\gl_m$-modules and
hence also of $\sll_m$-modules.

If $V \in \{P_m, P_m^*\}$, then as $\sll_m$-modules,
$\sll(m|1)_{\bar 1}\cong \C^m\oplus (\C^m)^*$ while 
$V_{\bar 1}$ is isomorphic to
the direct sum of $2m$ copies of the trivial module. 
Schur's lemma implies that $p|_{\sll(m|1)_{\bar 1}}=0$ and hence
$q=\widetilde{\ev}_a$ is an evaluation map of $\sll(m|1)$.

If $m = 3$ and $V = Q_3$, then as $\sll_3$-modules,
$\sll(3|1)_{\bar 1}\cong \C^3\oplus (\C^3)^*$ while 
$V_{\bar 1}$ is isomorphic to
the direct sum of $(\C^3)^*$ and $3$ copies of the trivial module. 
Moreover, the central element $E_{11} + E_{22} + E_{33}+ 3 E_{44}$ 
acts on the copy of $(\C^3)^*$ in $\sll(3|1)_{\bar 1}$ as multiplication 
by zero, while it acts on the copy of $(\C^3)^*$ in $(Q_3)_{\bar 1}$ 
as multiplication by 2. Applying Schur's lemma as above completes the 
argument in this case. The case when $m = 3$ and $V=(Q_3)^*$ is dealt with
in a similar manner.
\end{proof}

\begin{prop}\label{prop:4.3}
There are no bijective evaluation maps of $\sll(m|1)$ associated  with
the modules $P_{m}$, $P_{m}^*$, for $m \geq 3$, $Q_3$, and $Q_3^*$.
\end{prop}

\begin{proof} First we show that there are no bijective evaluation maps
associated with $P_m$ for $m \geq 3$.
Since $P_{m}=m\,\C^{m|1}\oplus m\, \ttr$, 
any point $a$ of $(P_{m})_{\0}$ is annihilated by the odd positive root
vector $E_{1,m+1}$, i.e., $\ev_a(E_{1,m+1})=E_{1,m+1}\cdot a=0$. Thus,
$\ev_a$ associated with the module  $P$ is not bijective for any point
$a\in (P_{m})_{\0}$. Similarly, the odd negative root vector $E_{m+1,1}$
annihilates any point of $(P_{m}^*)_{\0}$ and hence there are not
bijective evaluation maps associated with the module $P^*_m$.

Assume now that $m = 3$, $V = Q_3$, and $a \in V_{\bar 0}$.
A direct calculation shows that 
\[\dim \Span \{\ev_a(E_{14}),  \ev_a(E_{24}), \ev_a(E_{34})\} \leq 2 < 3 =
\dim \Span \{E_{14}, E_{24}, E_{34}\} \, ,\]
proving that $\ev_a$ is not a bijective evaluation map. The case when
$V = Q_3^*$ is dealt with in a similar way.
\end{proof}

\begin{proof}[Proof of Theorem \ref{thm1}]
Let $m\geq 3$. Assume to the contrary that $\LL$  is an LSSA 
on $\sll(m|1)$.
Let $V$ be the $m^2|2m$-dimensional $\sll(m|1)$-module given by $\LL$.
Then there exists
a bijective 1-cocycle of $\sll(m|1)$ associated with  $V$, and $V_{\0}$
induces  an LSA $\LL_{\0}$ on $\gl_m$.
It follows from \cite[Theorem 4.5]{Bau1999} that $V_{\0}$, as an
$\sll_m$-module, is isomorphic to $m\,\C^m$ or $m\,(\C^m)^*$.
Hence, by Proposition \ref{prop:4.4}, $V$ is  isomorphic to one 
of $P_{m},P_{m}^*,Q_3,$ or $Q_3^*$. 
Proposition \ref{prop:4.3} completes the proof.
\end{proof}

%%%%%%%%%%%%%%%%%%%%%%%%%%%%%%%%%%%%%%%%%%%SECTION 1
\section{Proof of Theorem \ref{thm2}} \label{sec5}
\setcounter{equation}{0}
\renewcommand{\theequation}
{5.\arabic{equation}}
\setcounter{theorem}{0}
\renewcommand{\thetheorem}
{5.\arabic{theorem}}

Throughout this section $\g = \sll(2|1)$ and we use the notation
introduced in Section \ref{sec3} some of which we recall for convenience.
Let $h=E_{11}-E_{22}$ and $z=E_{11}+E_{22}+2E_{33}$. A dominant integral
weight $\lambda$ of $\g$ is of the form $(i,k) \in \Z_{\geq 0} \times
\C$, where $i = \lambda(h)$ and $k = \lambda(z)$; the corresponding
irreducible highest weight module and Kac module are denoted respectively
by $V(i,k)$ and $K(i,k)$. The weight 
$(i,k)$ is atypical if and only if $k = i$ or $k=-i-2$. We index the
atypical irreducible $\g$-modules by $\Z$: $V_i := V(i,i)$ for $i \geq 0$
and $V_{i} = V(-i-1,  i-1)$ for $i < 0$. Finally,
$S_j$ denotes the $j+1$-dimensional irreducible $\sll_2$-module.

The proof of Theorem \ref{thm2} is carried out in the rest of this
section. Namely, in Section \ref{sec5.1} we construct the LSSAs $\A_k$,
$\B_{k_1,k_2}$, and $\CC_k$.
In Section \ref{sec5.2} we describe the $4|4$-dimensional modules 
which may be associated with LSSAs and prove that every LSSA on
$\sll(2|1)$ is isomorphic 
to an LSSA among $\A_k$, $\B_{k_1,k_2}$, and $\CC_k$. Finally, in Section
\ref{sec5.3} we establish the isomorphisms among $\A_k$, $\B_{k_1,k_2}$,
and $\CC_k$.

\subsection{The LSSAs $\A_k$, $\B_{k_1,k_2}$, and $\CC_k$} 
\label{sec5.1}
In this section we define the LSSAs $\A_k$, $\B_{k_1,k_2}$, and $\CC_k$
by providing a $\g$-module $M$ along with a vector $a \in M_{\0}$ for which
the evaluation map $\ev_a$ is bijective.

\subsubsection{{\bf $\A_k$ for $k \in \C \setminus \{-1, -3\}$}}
\label{sec5.1.1} Let $M = K(1, k)$. Then, as an $\sll_2$-module, 
$M_{\0} \cong S_1 \oplus S_1$. Moreover, 
as a $\gl_2$-module, $M_{\0} = M_1 \oplus M_2$, where $z$ acts on $M_1$
and $M_2$ as multiplication by $k$ and $k+2$ respectively. Denoting 
the highest weight vector of $M$ by $v_0$, we note that
\[\{v_0, v_1:= E_{21}v_0, w_0:=E_{31}E_{32}v_0, w_1:= E_{21}E_{31}E_{32}v_0\}\]
is a basis of $M_{\0}$. Similarly,
\[\{E_{32}v_0, E_{21}E_{32}v_0, E_{21}^2E_{32}v_0, E_{31}v_0\}\] 
is a basis of $M_{\w}$.

Consider $a := v_0 + w_1$. The
the action of $\g$ on $a$ is as follows:
\begin{equation} \label{eq510}
\begin{array}{ccc} 
h \cdot a & = & v_0 - w_1\\
z \cdot a & = & kv_0 + (k+2) w_1\\
E_{12} \cdot a & = & w_0\\
E_{21} \cdot a & = & v_1\\
E_{13} \cdot a & = & \frac{k+1}{2} E_{21}E_{32}v_0 +  \frac{k-1}{2}
E_{31}v_0\\
E_{23} \cdot a & = & \frac{k+3}{4} E_{21}^2E_{32}v_0\\
E_{31} \cdot a & = & E_{31}v_0 \phantom{\ .} \\
E_{32} \cdot a & = & E_{32}v_0 \ .\\
\end{array}
\end{equation}
It is immediate that, for $k \neq -1, -3$, the vectors in the right hand
side of \eqref{eq510} are linearly independent and hence $\ev_a$ is
bijective.
Thus, for $k \neq -1, -3$, the pair $(M,a)$ defines an LSSA on $\sll(2|1)$
which we denote by $\A_k$.

\subsubsection{{\bf $\B_{k_1,k_2}$ for $k_1, k_2 \in \C \setminus \{0\},$
$k_1 + k_2 \neq -2$}} \label{sec5.1.2}
In this case we set $M := \Pi K(0, k_1) \oplus \Pi K(0, k_2)$, where, as
usual $\Pi$ stands for the change-of-parity functor. As a $\gl_2$-module,
$M_{\0} = M_1 \oplus M_2$, where
$\dim M_1 = \dim M_2 = 2$ and $z$ acts on $M_1$ and $M_2$ as multiplication by
$k_1+1$ and $k_2+1$ respectively. As above, let $v_0$ be the highest
weight vector of $M_1$ and $w_1$ 
be the lowest weight vector of $M_2$.
Exactly as in the case of $\A_k$, one checks that, for $a = v_0 + w_1$,
the map $\ev_a$ is bijective as long as $k_1 + k_2 \neq -2$ and 
$k_1, k_2 \neq 0$. We leave completing the details to the reader. 
The resulting LSSA is denoted by $\B_{k_1, k_2}$.

\subsubsection{{\bf $\CC_{k}$ for $k \in \C \setminus \{0, -1\}$}}
\label{sec5.1.3}
In this case we set $M := \Pi K(0, k)^{(2)}$. As a $\gl_2$-module,
$M_{\0}$ is a non-trivial extension of $M_1$ by $M_2$, where 
$M_1 \cong M_2$, 
$\dim M_1 = \dim M_2 = 2$ and $z$ acts on $M_1$ and $M_2$ as
multiplication by $k+1$. Note that $z$ acts on the (two-dimensional)
highest weight space of $M_{\0}$ by the
$2 \times 2$-matrix $ \left(\begin{array} {cc} k+1 & 0\\ 1 & k+1
\end{array} \right)$. Let $v_0$ be a preimage in $M_{\0}$ of the highest
weight vector of $M_1$ and $w_1$ be the
lowest weight vector of $M_2 \subset M_{\0}$.
Exactly as in the case of $\A_k$, one checks that, for $a = v_0 + w_1$,
the map $\ev_a$ is bijective as long as and $k \neq 0, -1$. We leave completing
the details to the reader. The resulting LSSA is denoted by $\CC_{k}$.

%\begin{rem} \label{rem5.15}
%There is an LSSA $\B_{k_1,k_2}$ associated with the module 
%$M := \Pi K(0, k_1) \oplus \Pi K(0, k_2)$ even if one or both
%of the parameters $k_1$ and $k_2$ equal zero, as long as 
%$k_1 + k_2 \neq -2$. To construct these LSSAs one needs to
%use evaluation maps at different base points. Likewise,
%there is an LSSA $\CC_0$. In fact, in the Appendix we provide 
%explicit multiplication tables which are valid even for the parameters
%excluded in \ref{sec5.1.2} and \ref{sec5.1.3}.
%\end{rem}

\subsection{$\g$-modules associated with LSSAs} \label{sec5.2} Let $\LL$
be an LSSA on $\g = \sll(2|1)$ and let $M$ denote the corresponding
$\g$-module. Since $\LL_{\0}$ is a left-symmetric algebra with corresponding
$\gl_2$-module $M_{\0}$, Baues's classification theorem implies that, as
an $\sll_2$-module $M_{\0} \cong S_3$ or $M_{\0} \cong S_1 \oplus S_1$.
This fact, along with $\dim M = 4|4$,
imply that the composition factors of $M$ are among the following modules (cf.
Section \ref{sec3.4}):
\[K(1,k), k \neq 1, -3, \quad \Pi K(0,k), k \neq 0, -2, \quad \Pi V_{-3},
\quad V_{-2}, \quad \Pi V_{-1}, \quad \Pi V_0, \quad V_1,\quad \Pi V_2,
\quad V_3.\]
More precisely, to obtain a module of dimension $4|4$, we need to combine
composition factors from one of the following sets:
\begin{enumerate}
\item $\{K(1,k)\}$, $k \neq 1, -3$;
\item $\{\Pi V_{-3}, \Pi V_0\}$ or $\{V_3, \Pi V_0\}$;
\item $\{V_{-2}, \Pi V_{-1}\}$, $\{V_{-2}, V_{1}\}$, $\{\Pi V_{2}, \Pi
V_{-1}\}$, or $\{\Pi V_{2},  V_{1}\}$;
\item $\{\Pi K(0, k_1), \Pi K(0,k_2)\}$, $k_1, k_2 \not \in \{0, -2\}$;
\item $\{\Pi K(0, k), \Pi V_{-1}, \Pi V_0\}$ or $\{\Pi K(0, k), V_{1}, \Pi
V_0\}$; 
\item $\{\Pi V_{-1}, \Pi V_{-1}, \Pi V_0, \Pi V_0\}$, $\{\Pi V_{-1}, V_{1},
\Pi V_0, \Pi V_0\}$,  or $\{V_{1},  V_{1}, \Pi V_0, \Pi V_0\}$.
\end{enumerate}

First we prove Theorem \ref{thm2} (1):

\begin{prop}\label{lem5.2}
Let $\LL$ be an LSSA with corresponding $\g$-module $M$ and 1-cocycle
$q$. Then $q$ is an evaluation map.
\end{prop}

\begin{proof} The list of possible composition factors of $M$ above shows
that every composition factor of $M_{\w}$ considered as an
$\sll_2$-module is isomorphic to $S_0$ or $S_2$. 
Noting that $\sll(2|1)_{\bar 1}$, considered as an $\sll_2$-module, is
isomorphic to $S_1 \oplus S_1$, an argument as in Lemma \ref{lemma:4.2}
proves that $q$ is an evaluation map.
\end{proof}

Proposition \ref{lem5.2} and Lemma \ref{lem:3.5} imply immediately:

\begin{coro} \label{cor510}
If $\LL_1$ and $\LL_2$ are two LSSAs on $\g$ corresponding to the same
$\g$-module $M$, then $\LL_1 \cong \LL_2$.
\end{coro}

\begin{rem} \label{rem511}
Recall from Section \ref{sec3.1} 
that $M^{-st}$ is isomorphic to the twist of the $\g$-module $M$ by the outer
automorphism $-st$ of $\g$. If $(M, q)$ is the pair of a $\g$-module and
a bijective 1-cocycle corresponding to
an LSSA $\LL$, then the LSSA $\LL^{st}$ corresponding to the pair
$(M^{-st}, q \circ (-st))$ is isomorphic to $\LL$.
In particular, to list all LSSAs it suffices to determine which
$\g$-modules $M$ (up to a twist by $-st$) admit  bijective evaluation
maps.
\end{rem}

\begin{lem}\label{lem5.3}
Let $M$ be a $\g$-module whose composition factors are in the list above.
Assume $M$ satisfies one of the conditions:
\begin{enumerate}
\item  Both $\Pi V_{-1}$ and $V_1$ are composition factors of $M$;
\item  $\Pi V_0$ is a quotient of $M$;
\item  $\Pi V_{-1}$ or $V_1$ is a submodule of $M$.
\end{enumerate}
Then there is no bijective evaluation map associated with $M$.
\end{lem}

\begin{proof}
(1) If both $\Pi V_{-1}$ and $V_1$ are composition factors of $M$, then,
as a $\gl_2$-module, $M_{\0} = M_1 \oplus M_2$, where both $M_1$ and
$M_2$ are 2-dimensional irreducible
$\sll_2$-modules and $z$ acts on $M_1$ as multiplication by $1$ and on
$M_2$ -- by $-1$. Hence, $M_2 \cong M_1^*$, which implies that the
$\gl_2$-module 
$M_{\0}$ does not admit a bijective evaluation map.

(2) Assume $M'$ is a submodule of $M$ such that $M/M' \cong \Pi V_0$.
Then the image of $\ev_a$ is contained in $M'$ and hence $\ev_a$ is not
bijective.

(3) Assume that $V_1$ is a submodule of $M$ and consider the list of
possible composition
factors of $M$ above. If the composition factors of $M$ are 
$\{V_{-2}, V_{1}\}$, then, for 
any $a \in M_{\0}$,  $\ev_a(E_{13}) = \ev_a(E_{23}) = 0$  and thus
$\ev_a$ is not bijective. If these are $\{\Pi V_{2},  V_{1}\}$, then, for 
any $a \in M_{\0}$,  $\ev_a(E_{31})$ and $\ev_a(E_{32})$ are  linearly
dependend and thus $\ev_a$ is not bijective. 
If these are $\{\Pi K(0, k), V_{1}, \Pi V_0\}$ or $\{V_{1},  V_{1}, \Pi V_0,
\Pi V_0\}$, then $\Pi V_0$ is a quotient of $M$ and we refer to (2). 
Finally, if these are $\{\Pi V_{-1}, V_{1}, \Pi V_0, \Pi V_0\}$, then we refer to
(1). The case when $\Pi V_{-1}$ is a submodule of $M$ is dealt with in a
similar way.
\end{proof}

We are now ready to prove Theorem \ref{thm2} (2).

\begin{prop} \label{prop5.15} Let $\LL$ be an LSSA. Then $\LL$ is
isomorphic to one of 
$\A_{k},  k\in\C\setminus\{-1,-3\}$, $\B_{k_1,k_2},   k_1 + k_2 \neq -2,
k_1, k_2 \in \C\setminus\{0\}$, or  $\CC_{k},  
k\in\C\setminus\{-1,0\}$.
\end{prop}

\begin{proof}
Let $\LL$ be an LSSA with corresponding $\g$-module $M$. In view of
Corollary \ref{cor510} and Remark \ref{rem511}, it suffices to show that
$M$ or $M^{-st}$ 
is isomorphic to a module corresponding to one of the LSSAs $\A_{k}, 
k\in\C\setminus\{-1,-3\}$, $\B_{k_1,k_2},   k_1 + k_2 \neq -2, k_1, k_2
\in \C\setminus\{0\}$, and  
$\CC_{k},   k\in\C\setminus\{-1, 0\}$. We consider the six cases listed
above for the composition factors of $M$.

(1) $M = K(1,k)$ for $k \neq 1, -3$ is irreducible. Then, for $k \neq
-1$, $M$ gives rise to the LSSA $\A_k$. Consider $M = K(1,-1)$. In this
case, as a $\gl_2$-module, 
$M_{\0} = M_1 \oplus M_2$, where $M_2 \cong M_1^*$.  This means that
$M_{\0}$ does not admit a bijective evaluation map and hence $M=K(1,-1)$
does not give rise to an LSSA.

(2) Since there are no extensions between $\Pi V_0$ and $V_3$ or $\Pi
V_{-3}$, in this case $\Pi V_0$ is a quotient of $M$ and, by Lemma
\ref{lem5.3}, $M$ does not give rise to an LSSA.

(3) Since there are no extensions between $V_1$ and $V_{-2}$ or between
$\Pi V_{-1}$ and $\Pi V_2$, Lemma \ref{lem5.3} excludes these cases. For
the other two pairs, again
Lemma \ref{lem5.3}, leaves only two possible modules: a non-trivial
extension of $V_1$ by $\Pi V_2$ and a non-trivial extension of $\Pi
V_{-1}$ by $V_{-2}$. In the former case
$M \cong K(1, 1)$ and hence corresponds to the LSSA $\A_1$. In the latter
case $M^{-st} \cong K(1,1)$ and, by Remark \ref{rem511}, it corresponds to 
an LSSA isomorphic to $\A_1$.

(4) In this case we need to consider two cases: when $M$ is completely
reducible and when $M$ is indecomposable.

Assume first that $M = \Pi K(0, k_1) \oplus \Pi K(0, k_2)$, where $k_1,
k_2 \not \in \{0, -2\}$. If $k_1 + k_2 \neq -2$, $M$ corresponds to the
LSSA $\B_{k_1, k_2}$. 
If $k_1 + k_2 = -2$, then $M_{\0}$ is a self-dual $\gl_2$-module and hence
it does not admit a bijective evaluation map. In particular, for $k_1 +
k_2 = -2$, 
$M = \Pi K(0, k_1) \oplus \Pi K(0, k_2)$ does not give rise to an LSSA.

Next assume that $M$ is indecomposable. Then $k_1 = k_2 =: k \neq 0, -2$
and $M \cong \Pi K(0,k)^{(2)}$. Thus, for $k \neq -1$, $M$ corresponds to
the LSSA $\CC_{k}$. 
For $k = -1$, the $\gl_2$-module $M_{\0}$ is self-dual; hence it does not
admit a bijective evaluation map and nor does $M$.

(5) If the composition factors of $M$ are $\{\Pi K(0, k), \Pi V_{-1}, \Pi
V_0\}$ with $k \neq 0, -2$, Lemma \ref{lem5.3} implies that $M = M_1 \oplus M_2$, where $M_1
\cong \Pi K(0, k)$ and 
$M_2 \cong \Pi K(0,-2)$. Thus $M$ gives rise to $\B_{k, -2}$. If the
composition factors of $M$ are $\{\Pi K(0, k), V_{1}, \Pi V_0\}$ with $k \neq 0, -2$, then the
composition factors of $M^{-st}$ are
$\{\Pi K(0, -k-2), \Pi V_{-1}, \Pi V_0\}$ and, as above,
$M^{-st}$ (and hence $M$) corresponds to an LSSA isomorphic to 
$\B_{-k-2, -2}$.

(6) Assume that the composition factors of $M$ are 
$\{\Pi V_{-1},  \Pi V_{-1}, \Pi V_0, \Pi V_0\}$.  Lemma \ref{lem5.3} 
implies that $M$ is isomorphic to $\Pi K(0, -2) \oplus \Pi K(0, -2)$ or to
$\Pi K(0, -2)^{(2)}$. The former module gives rise to $\B_{-2,-2}$,
while that latter gives rise to $\CC_{-2}$.
If the composition factors of $M$ are $\{V_{1},  V_{1}, \Pi V_0, \Pi V_0\}$, 
then $M^{-st}$ is isomorphic to $\Pi K(0, -2) \oplus \Pi K(0, -2)$ or to
$\Pi K(0, -2)^{(2)}$ and hence $M$ gives rise to an LSSA isomorphic to 
$\B_{-2,-2}$ or $\CC_{-2}$.
Finally, Lemma \ref{lem5.3} shows that a module with composition factors
$\{\Pi V_{-1},  V_{1}, \Pi V_0, \Pi V_0\}$ does not give rise to an LSSA.
\end{proof}

\subsection{Isomorphisms} \label{sec5.3}
We complete this section with the proof of Theorem \ref{thm2} (3).

\begin{prop} \label{prop530}
$\A_k \cong \A_{-2-k}$, $\B_{k_1,k_2} \cong \B_{k_2,k_1} \cong
\B_{{-2-k_1},{-2-k_2}}\cong \B_{{-2-k_2},{-2- k_1}}$, and $\CC_k \cong
\CC_{-2-k}$.
Moreover, these are the only isomorphisms among $\A_k$,
$\B_{k_1,k_2}$, and $\CC_k$.
\end{prop}

\begin{proof} Let $\LL$ and $\LL'$ be two LSSAs with corresponding
modules $M$ and $M'$ respectively. Remark \ref{rem511} implies that 
$\LL' \cong \LL$ if and only if $M' \cong M$
or $M' \cong M^{-st}$. Applying Proposition \ref{prop3.11} to the list of
$\g$-modules corresponding to the LSSAs $\A_k, \B_{k_1,k_2}, \CC_k$
 completes the proof.
\end{proof}

%\begin{rem} \label{rem5.33}
%The apparent inconsistency between the isomorphism $\A_k \cong \A_{-2-k}$
%and the fact that $\A_k$ is defined for $k = 1$ but not for $k=-3$ is due 
%to the fact that the $\g$-module $M = K(1,1)$ which corresponds to $k=1$
%is not irreducible and $M^{-st}$ is not a Kac module. Moreover, 
%as the proof of Theorem \ref{thm2} shows, 
%the module $K(1,-3)$ does not give rise to an LSSA structure.
%\end{rem}

%%%%%%%%%%%%%%%%%%%%%%%%%%%%%%%%%%%%%%%%%%%SECTION 1

\section{An Example of LSSAs on $\sll(m+1|m)$}
\setcounter{equation}{0}
\renewcommand{\theequation}
{6.\arabic{equation}}
\setcounter{theorem}{0}
\renewcommand{\thetheorem}
{6.\arabic{theorem}}

Throughout this section $m$ is a fixed positive integer.
Before we prove Theorem \ref{thm4} by providing an 
$\sll(m+1|m)$-module $U$ which admits a bijective evaluation map, 
we recall some facts about the exterior square of a super vector space.

Let $W = W_{\0} \oplus W_{\w}$ be a $\Z_2$-graded vector space.
Then the exterior square of $W$ is by definition
\[\wedge^2 W := (W \otimes W)/\Span\{u\otimes v + (-1)^{|u| \, |v|} v \otimes u\}.\]
As non-graded vector spaces, $(\wedge^2 W)_{\0} = \wedge^2 W_{\0} \oplus S^2 W_{\w}$
and $(\wedge^2 W)_{\w} = W_{\0} \otimes W_{\w}$, where $S^2 W_{\w}$ is the usual
symmetric square of $W_{\w}$. In particular, 
\[\dim\, \wedge^2 W = \binom{\dim W_{\0}}{2} + \binom{\dim W_{\w}+1}{2} \vert
\dim W_{\0} \dim W_{\w}\ .\]
If $u, v \in W$, we denote the image of $u \otimes v$ in $\wedge^2 W$ under the 
natural projection $W \otimes W \twoheadrightarrow \wedge^2 W$ by $uv$.
If 
$\{e_1, e_{2}, \dots, e_p\}$ is a basis of $W_{\bar 0}$ and
$\{\xi_1, \xi_{2}, \dots, \xi_{q}\}$ is a basis of $W_{\bar 1}$, then 
$\{e_i e_j, \xi_s\xi_t \, | \, 1 \leq i < j \leq p, 1 \leq s \leq t \leq q\}$
is a basis of $(\wedge^2 W)_{\0}$ and
$\{e_i \xi_s \, | \, 1 \leq i \leq p, 1 \leq s \leq q\}$ is a basis of 
$(\wedge^2 W)_{\w}$.

\begin{proof}[Proof of Theorem  \ref{thm4}]
Let $W := \C^{m+1|m}$ be the standard module of $\sll(m+1|m)$. Set
$$U:=\Pi(\wedge^2 W)\oplus \Pi(\wedge^2 W), $$
i.e., $U$ is the direct sum of two copies of the exterior square
of $W$ with the parity reversed.

Let $\{e_1, e_{2}, \dots, e_{m+1}\}$  and
$\{\xi_1, \xi_{2}, \dots, \xi_{m}\}$ denote the standard bases of $W_{\0} = \C^{m+1|m}_{\bar 0}$ and 
$W_{\w} = \C^{m+1|m}_{\bar 1}$
respectively.
Given an element $v\in \Pi(\wedge^2 W)$, we denote  by $v'$ and
$v''$ the elements $(v,0)$ and $(0,v)$ in $U$. In this notation
\[
\{e'_i \xi'_s,\, \,  e''_i \xi''_s \, | \, 1 \leq i \leq m+1, 1 \leq s \leq m\}
\]
and
\[
\{e'_i e'_j, e''_i e''_j, \xi'_s \xi'_t, \xi''_s \xi''_t \, | \, 1 \leq i < j \leq m+1,
1 \leq s \leq t \leq m\}
\]
are bases of $U_{\0}$ and $U_{\w}$ respectively.

Consider the element $$a:=
\sum_{i=1}^{m} (e'_{i+1}\xi'_{i}+e''_{i}\xi''_{i}) \in U_{\bar 0}.$$
We show below that $\ev_{a}:\sll(m+1|m)\ra U$ is a bijective evaluation map,
thus proving Theorem \ref{thm4}.
Since $\dim \sll(m+1|m))=\dim U$, to prove that $\ev_{a}$ is
bijective, it suffices to prove it is injective, i.e, that $\ker \ev_a = 0$.

Assume $X = \left( \begin{array}{cc} A & B\\ C & D \end{array}\right) \in \sll(m+1|m)$,
where $A,B,C,$ and $D$ are matrices of sizes 
$(m+1) \times (m+1), (m+1) \times m, m \times(m+1),$ and $m \times m$ respectively.
Let $A = (a_{i,j}), B= (b_{i,s}), C = (c_{s,i})$, and 
$D = (d_{s,t})$ with $1 \leq i,j \leq m+1$
and $1 \leq s,t \leq m$. We calculate
\[
\begin{array}{rl}
\ev_a(X)& = \sum_{i=1}^{m+1} \sum_{s=1}^{m} a_{i,s+1}e'_i\xi'_s + 
\sum_{i=1}^{m} \sum_{s=1}^{m} d_{s,i} \, e'_{i+1}\xi'_s \\
&\\
& + \sum_{i=1}^{m+1} \sum_{s=1}^{m} a_{i,s}e''_i\xi''_s + 
\sum_{i=1}^{m} \sum_{s=1}^{m} d_{s,i} \, e''_i \xi''_s \\
&\\
& + \sum_{i=1}^{m} \sum_{j=1}^{m+1} b_{j,i} \, e'_{i+1} e'_j +
\sum_{i=1}^{m} \sum_{j=1}^{m+1} b_{j,i} \, e''_i e''_j \\
& \\
&+ \sum_{s=1}^{m} \sum_{t=1}^{m} c_{s,t+1} \xi'_s \xi'_t +
\sum_{s=1}^{m} \sum_{t=1}^{m} c_{s,t} \, \xi''_s \xi''_t \ .
\end{array}
\]
The expression above shows that $X \in \ker \ev_a$ if and only if
\[
\begin{array}{lllcl}
a_{1, s+1} = 0 & & a_{i+1, s+1} + d_{s,i} = 0& {\text{for}} & 1 \leq i,s \leq m,\\
&&&&\\
a_{m+1, s} = 0 & & a_{i,s} + d_{s,i} = 0& {\text{for}} & 1 \leq i,s \leq m, \\
&&&&\\
b_{1,j} = 0&& b_{i+1,j} - b_{j+1,i} = 0& {\text{for}} & 1 \leq i \neq j \leq m,\\
&&&&\\
b_{m+1,j} = 0&& b_{i,j} - b_{j,i} = 0& {\text{for}} & 1 \leq i \neq j \leq m,\\
&&&&\\
c_{s,s+1} = 0 && c_{s,t+1} + c_{t,s+1} = 0 & {\text{for}} & 1 \leq s \neq t \leq m,\\
&&&&\\
c_{s,s} = 0 && c_{s,t} + c_{t,s} = 0 & {\text{for}} & 1 \leq s \neq t \leq m.\\
\end{array}
\]
An easy and somewhat tedious calculation shows that the solutions of the system above
are the matrices $A=c I_{m+1}, D = - c I_m$, $B= 0$, $C=0$, where $c$ is a scalar.
Since $\text{str}\, X = (2m+1) c$ and $X \in \sll(m+1|m)$, 
we conclude that $c = 0$, i.e., $X=0$. 
This proves that $\ev_a$ is  injective and completes the proof of the theorem.
\end{proof}

%\clearpage

\section*{Appendix}
Set
\begin{eqnarray*}
&&x_{1}:=E_{12}, \ x_{2}:=E_{21}, \ x_{3}:=E_{11}-E_{22}, \ x_{4}:=
E_{11}+E_{22}+2E_{33},\\
&&y_{1}:=E_{31}, \ y_{2}:=E_{32}, \ y_{3}:=E_{13}, \ y_{4}:=E_{23}.
\end{eqnarray*}
Below we provide the multiplication tables for $\A_k, \B_{k_1,k_2}$, and $\CC_k$.
%\vspace{4em}

\begin{table}[!ht]
 \Small
\caption{The LSSAs $\A_{k},~ k\neq -1,-3$  } % title name of the table
\centering % centering table
\begin{tabular}{|c|c|c|c|c|c|c|c|c|}
\hline
&   $x_1$   &     $x_2$      &   $x_3$   &     $x_4$    \\
 \hline
 $x_1$   &   0  &     $\frac{1}{2(k+1)}((k+2)x_3+x_4)$      &   $-x_1$   
&     $(k+2)x_1$    \\
 $x_2$&   $\frac{1}{2(k+1)}(-kx_3+x_4)$    &     0   &   $x_2$  
&  $kx_2$   \\
 $x_3$   &   $x_1$   &     $-x_2$      &   $\frac{1}{k+1}(x_3+x_4)$   
&     $\frac{1}{k+1}\big(k(k+2)x_3-x_4\big)$        \\
 $x_4$   &   $(k+2)x_1$   &    $kx_2$   &
$\frac{1}{k+1}\big(k(k+2)x_3-x_4\big)$  &    
$\frac{-k(k+2)}{k+1}x_3+\frac{k^2+2k+2}{k+1}x_4$   \\
$y_1$ &    $0$&     $\frac{4}{k+3}y_4$&     $y_1$&    $ky_1$      \\
$y_2$  &    $0$&     $\frac{2}{k+1}(y_1-2y_3)$&     $y_2$&    $ky_2$   \\
$y_3$  &    $-\frac{k+3}{4}y_2$&     $0$&     $-y_3$&    $(k+2)y_3$ \\
$y_4$  & $\frac{(k+3)(k-1)}{4(k+1)}y_1+\frac{2}{k+1}y_3$ &  0   &    
$-y_4$&    $(k+2)y_4$ \\
\hline
&   $y_1$&     $y_2$&     $y_3$&    $y_4$   \\
 \hline
$x_1$     &    $-y_2$&     0&     $-\frac{k+3}{4}y_2$&  
$\frac{k+3}{k+1}(\frac{k-1}{4}y_1+y_3)$   \\
$x_2$ &    $\frac{4}{k+3}y_4$ &  -$\frac{1}{k+1}\big((k-1)y_1+4y_3\big)$
&   $y_4$&  0  \\
$x_3$   &    0&     $2y_2$&     0&   $-2y_4$      \\
$x_4$    &   $(k+1)y_1$  &     $(k+1)y_2$  &     $(k+1)y_3$&   
$(k+1)y_4$   \\
$y_1$  &   $0$   &     $-2x_1$      &   $\frac{1}{4}(x_4-kx_3)$   & 0 \\
$y_2$  &  $2x_1$    &     $0$      &   $x_1$ &  
$\frac{k+3}{4(k+1)}(x_4-kx_3)$  \\
$y_3$  & $\frac{1}{4}\big((k+2)x_3+x_4\big)$ &    0      &   $0$   &   
$\frac{(k+3)(k-1)}{8}x_2$  \\
$y_4$  & $x_2$  &     $\frac{k-1}{4(k+1)}\big((k+2)x_3+x_4\big)$    &
$-\frac{(k+3)(k-1)}{8}x_2$   &     $0$\\
\hline
\end{tabular}
\end{table}
\vspace{4em}

\begin{table}[!ht]
 \tiny
\caption{The LSSAs $\B_{k_1,k_2},~ k_1+k_2\neq -2$} 
% title name of the table
\centering % centering table
\begin{tabular}{|c|c|c|c|c|c|c|c|c}
\hline
&   $x_1$   &     $x_2$      &   $x_3$   &     $x_4$    \\
 \hline
$x_1$   &   0  &     $\frac{1}{k_1+k_2+2}\big((k_2+1)x_3+x_4\big)$      
&   $-x_1$   &     $(k_2+1)x_1$    \\
$x_2$  &   $\frac{1}{k_1+k_2+2}\big(-(k_1+1)x_3+x_4\big)$    &     0   
&   $x_2$  &  $(k_1+1)x_2$   \\
$x_3$   &   $x_1$   &     $-x_2$      
&   $\frac{1}{k_1+k_2+2}\big((k_2-k_1)x_3+2x_4\big)$   &    
$\frac{2(k_1+1)(k_2+1)}{k_1+k_2+2}x_3+\frac{k_1-k_2}{k_1+k_2+2}x_4$ \\
$x_4$   &   $(k_2+1)x_1$   &    $(k_1+1)x_2$   
&   $\frac{2(k_1+1)(k_2+1)}{k_1+k_2+2}x_3+\frac{k_1-k_2}{k_1+k_2+2}x_4$ 
&     $\frac{(k_1+1)(k_2+1)(k_1-k_2)}{k_1+k_2+2}x_3+
\frac{(k_1+1)^2+(k_2+1)^2}{k_1+k_2+2}x_4$   \\
$y_1$  &    $y_2$&     0&     $y_1$&    $(k_1+1)y_1$      \\
$y_2$  &    $0$&     $y_1$&     $-y_2$&    $(k_2+1)y_2$   \\
$y_3$  &    $0$&     $-y_4$&     $-y_3$&    $(k_2+1)y_3$ \\
$y_4$  & $-y_3$ &  0   &     $y_4$&    $(k_1+1)y_4$ \\
\hline
&   $y_1$&     $y_2$&     $y_3$&    $y_4$   \\
 \hline
$x_1$     &    0&     0&     0&    0   \\
$x_2$  &    0&     0&   0&  0  \\
$x_3$     &    0&     0&     0&    0      \\
$x_4$    &   $(k_1+2)y_1$  &     $(k_2+2)y_2$  &     $k_2y_3$& 
$k_1y_4$   \\
$y_1$  &   $0$   &     0      &  
$\frac{k_2}{2(k_1+k_2+2)}\big(x_4-(k_1+1)x_3\big)$   &  
$-\frac{k_1}{2}x_2$ \\
$y_2$   &   0   &     $0$      &   $-\frac{k_2}{2}x_1$ &  
$\frac{k_1}{2(k_1+k_2+2)}\big((k_2+1)x_3+x_4\big)$ \\
$y_3$   &  $\frac{k_1+2}{2(k_1+k_2+2)}\big((k_2+1)x_3+x_4\big)$  &  
$\frac{k_2+2}{2}x_1$    &   $0$   &     $0$  \\
$y_4$   &  $\frac{k_1+2}{2}x_2$  &   
$\frac{k_2+2}{2(k_1+k_2+2)}\big(x_4-(k_1+1)x_3\big)$ & $0$   &     $0$\\
\hline
\end{tabular}
$\phantom{ldldld}$\\
$\phantom{ldldld}$
\end{table}
\vspace{4em}

 \begin{table}[!ht]
\tiny
\caption{The LSSAs $\CC_{k},~ k\neq -1$ } % title name of the table
\centering % centering table
\begin{tabular}{|c|c|c|c|c|c|c|c|c|}
\hline
&   $x_1$   &     $x_2$      &   $x_3$   &     $x_4$    \\
 \hline
$x_1$   &   0  &      $\frac{1}{2(k+1)}(x_4-x_1)+\frac{1}{2}x_3$      
&   $-x_1$   &     $(k+1)x_1$    \\
$x_2$  &    $\frac{1}{2(k+1)}(x_4-x_1)-\frac{1}{2}x_3$    &     0   
&   $x_2$  &   $\frac{1}{2(k+1)}(x_4-x_1)+(k+1)x_2-\frac{1}{2}x_3$   \\
$x_3$   &   $x_1$   &     $-x_2$      &   $\frac{1}{k+1}(x_4-x_1)$   
&     $x_1+(k+1)x_3$        \\
$x_4$   &   $(k+1)x_1$   
&      $\frac{1}{2(k+1)}(x_4-x_1)+(k+1)x_2-\frac{1}{2}x_3$  
&   $x_1+(k+1)x_3$  &     $(k+1)(x_1+x_4)$   \\
$y_1$  &    $y_2$&     0&     $y_1$&    $(k+1)y_1+y_2$      \\
$y_2$  &    $0$&     $y_1$&     $-y_2$&    $(k+1)y_2$   \\
$y_3$  &    $0$&     $-y_4$&     $-y_3$&    $(k+1)y_3$ \\
$y_4$  & $-y_3$ &  0   &     $y_4$&    $(k+1)y_4-y_3$ \\
\hline
&   $y_1$&     $y_2$&     $y_3$&    $y_4$   \\
 \hline
$x_1$     &    0&     0&     0&    0   \\
$x_2$  &    0&     0&   0&  0  \\
$x_3$     &    0&     0&     0&    0      \\
$x_4$    &   $(k+2)y_1+y_2$  &     $(k+2)y_2$  &     $ky_3$&  
$ky_4-y_3$   \\
$y_1$  &   $0$   &     0&   $\frac{k}{4(k+1)}(x_4-x_1)-\frac{k}{4}x_3$  
&     $\frac{1}{4(k+1)}(x_1-x_4)-\frac{k}{2}x_2+\frac{1}{4}x_3$ \\
$y_2$   &   0   &     $0$      &   $-\frac{k}{2}x_1$ &    
$\frac{k+2}{4(k+1)}x_1+\frac{k}{4}x_3+\frac{k}{4(k+1)}x_4$  \\
$y_3$   & $\frac{k}{4(k+1)}x_1+\frac{k+2}{4}x_3+\frac{k+2}{4(k+1)}x_4$  
&     $\frac{k+2}{2}x_1$      &   $0$   &     $0$  \\
$y_4$   & $\frac{1}{4(k+1)}(x_4-x_1)+\frac{k+2}{2}x_2-\frac{1}{4}x_3$   &
$\frac{k+2}{4(k+1)}(x_4-x_1)-\frac{k+2}{4}x_3$      &   $0$   &     $0$\\
\hline
\end{tabular}
$\phantom{ldldld}$
 \end{table}
\vspace{3em}

Note that the tables above include the LSSAs $\B_{k, 0}$ with $k \neq -2$ (including $\B_{0,0}$) 
and $\CC_0$. These correspond to the $\sll(2|1)$-modules $(\Pi K(0,-2-k) \oplus \Pi K(0,-2))^{-st}$ and
$(\Pi K(0,-2)^{(2)})^{-st}$ respectively. In particular, they are isomorphic to $\B_{-2-k, -2}$ and 
$\CC_{-2}$ respectively.

%\section*{Declarations}
%\subsection*{Funding} This work was supported by NSFC (No. 11301061) and an NSERC Discovery Grant.
%\subsection*{Competing interests} The authors have no relevant financial or non-financial interests to disclose.
%\subsection*{Data availability statement} %The datasets generated  and/or analysed during the current study are available from the corresponding author on reasonable request. Data sharing is not applicable to this work as no datasets were generated or analyzed during the current study.

%\clearpage 

\end{document}